\documentclass[oneside]{amsart}
\newcommand{\formatswitch}{preprint}

\usepackage{amssymb}
\usepackage{ifthen}
\usepackage{verbatim}
\usepackage{psfig}
\usepackage[all]{xy}

\newcommand{\tref}[1]{(\ref{#1})}

\ifthenelse{\equal{\formatswitch}{big}}{
\setlength{\textwidth}{432pt}
\setlength{\oddsidemargin}{18.8775pt}

\setcounter{tocdepth}{1}
}{
\ifthenelse{\equal{\formatswitch}{paper}}{

\setcounter{tocdepth}{1}
}{
\setcounter{tocdepth}{2}
}
}
\DeclareMathAlphabet\EuScript{U}{eus}{m}{n}
\DeclareMathAlphabet\EuScriptb{U}{eus}{b}{n}

\newcommand{\sscr}[1]{\EuScript{#1}}

\newcommand{\claimenum}{\renewcommand{\theenumi}{\alph{enumi}}
 \renewcommand{\labelenumi}{\textit{(\theenumi)}}
 \renewcommand{\theenumii}{\roman{enumii}}
 \renewcommand{\labelenumii}{\textit{(\theenumii)}}
 \begin{enumerate}}
\newcommand{\claimenumend}{\end{enumerate}}

\newcommand{\romanenum}{\renewcommand{\theenumi}{\roman{enumi}}
 \renewcommand{\labelenumi}{\textit{(\theenumi)}}
 \renewcommand{\theenumii}{\alph{enumii}}
 \renewcommand{\labelenumii}{\textit{(\theenumii)}}
 \begin{enumerate}}
\newcommand{\romanenumend}{\end{enumerate}}

\newtheorem{dummy}{realdumb}[section]
\newtheorem{thm}{Theorem}

\newtheorem{lemma}[dummy]{Lemma}
\newtheorem{prop}[dummy]{Proposition}
{\theoremstyle{definition} }

{\theoremstyle{definition} }
\newtheorem{cor}{Corollary}[dummy]

\renewcommand{\text}{\mathrm}

\newcommand{\strutdepth}{\dp\strutbox}
\newcommand{\marginalnote}[1]
   {\strut\vadjust{\kern-\strutdepth\domarginalnote{#1}}}
\newcommand{\domarginalnote}[1]{\vtop to \strutdepth{
  \baselineskip\strutdepth
   \vss\llap{ #1\ \ }\null}}  
\newcounter{showlabelflag}
\newcounter{makelabelflag}
\newcommand{\showlabels}{\setcounter{showlabelflag}{1}}

\newcommand{\makelabels}{\setcounter{makelabelflag}{1}}
\newcommand{\hidelabels}{\setcounter{showlabelflag}{2}}
\newcommand{\mylabel}[1]{
  \ifthenelse{\value{makelabelflag}=1}
    {\label{#1}}{}
  \ifthenelse{\value{showlabelflag}=1}
    {\marginpar{#1}}{}\relax}

\newcommand{\R}{{\mathbf R}}

\newcommand{\N}{{\mathbf N}}

\ifthenelse{\equal{\formatswitch}{draft}}{
\showlabels
}{\relax}

\newcommand{\ZS}{Zappa-Sz\'ep }

\newcommand{\scr}{\sscr}

\newcommand{\zapprod}{\mathbin{\bowtie}}
\newcommand{\FB}{\scr{F}\zapprod B_\infty}
\newcommand{\FP}{\scr{F}\zapprod S_\infty}
\newcommand{\HB}{\scr{H}\zapprod B_\infty}
\newcommand{\HP}{\scr{H}\zapprod S_\infty}
\newcommand{\BZ}{B_\infty\zapprod\scr{Z}}
\newcommand{\PZ}{S_\infty\zapprod\scr{Z}}

\newcommand{\mymargin}[1]{
  \ifthenelse{\value{showlabelflag}=1}
    {\marginpar{#1}}{}\relax}

\newcounter{enumo}
\newcommand{\LeftAct}[2]{#1 \cdot #2}

\newcommand{\RRsh}{\kern -1 pt \Rsh}

\newcounter{keepitemnum}

\newcounter{keepitemnumm}

\begin{document}

\bibliographystyle{amsplain}
\begin{center}{\bfseries The Algebra of Strand Splitting. I.
\\
A Braided Version of Thompson's Group \protect\(V\protect\)\footnote{AMS
Classification (2000): primary 20F36, secondary 20B07,
20E99}}\end{center}
\vspace{3pt}
\begin{center}{MATTHEW G. BRIN}\end{center}
\vspace{4pt}
\vspace{3pt}
\begin{center}May 28, 2004\end{center}

\CompileMatrices


\makelabels
\hidelabels


\section{Introduction}\mylabel{IntroSec}

We construct a braided version \(BV\) of Thompson's group \(V\) that
surjects onto \(V\).  The group \(V\) is the third of three well
known groups \(F\), \(T\) and \(V\) created by Thompson in the 1960s
that have been heavily studied since.  See \cite{CFP} and Section 4
of \cite{brown:finiteprop} for an introduction to Thompson's groups.

The group \(V\) is a subgroup of the homeomorphism group of the
Cantor set \(C\).  It is generated by involutions \cite[Section
12]{brin:hd3} and, if the metric on \(C\) is ignored, \(V\) can be
viewed somewhat as a ``Coxeter group'' of permutations of \(C\).  In
\cite{brin:bv3} we find presentations for \(BV\) and \(V\) that
differ only in that the presentation for \(V\) has relations of the
form \(x^2=1\) that are not present in the presentation for \(BV\).
Thus \(BV\) can be thought of as an ``Artinification'' of \(V\).

Our motivation for creating \(BV\) is a relationship between \(BV\)
and the Thompson's groups \(F\) and \(V\) on the one hand, and
categories with multiplication on the other.  Given a category
\(\scr{C}\) with multiplication, an isomorphism expressing
associativity up to equivalence, and perhaps an isomorphism
expressing commutativity up to equivalence, there are groups and
epimorphisms \[i:G_1(\scr{C}) \rightarrow F,\quad\qquad
j:G_2(\scr{C}) \rightarrow V,\quad\qquad k:G_3(\scr{C}) \rightarrow
BV,\] that can be calculated from  \(\scr{C}\) and its attached
data for which \(i\) is an isomorphism if and only if \(\scr{C}\)
satisfies the axioms of a monoidal category, \(j\) is an isomorphism
if and only if \(\scr{C}\) satisfies the axioms of a symmetric,
monoidal category, and \(k\) is an isomorphism if and only if
\(\scr{C}\) satisfies the axioms of a braided tensor category.  See
\cite{joyal+street} and \cite{MacLane:categories} for definitions.
These results will be written up elsewhere.

As an intermediate step in understanding the group \(BV\), we also
construct a ``larger'' group \(\widehat{BV}\) that contains (and can
also be shown to be contained in) \(BV\) that is somewhat more
tractable.  If \(BV\) is regarded as a braided version of \(V\),
then \(\widehat{BV}\) is a braided version of \(\widehat{V}\) that
contains (and can also be shown to be contained in) \(V\) and that
is also somewhat more tractible than \(V\).

The group \(\widehat{V}\) acts on countably many copies of the
Cantor set, and another view of \(\widehat{BV}\) is that it is the
group \(B_\infty\) of finitary braids on countably many strands that
has been modified by allowing the strands to split and recombine.
This explains the first part of the title of this paper.  The group
\(BV\) is the subgroup in which all splitting, braiding and
recombining is confined to the first strand.  The group \(BV\) is
thus the ``braid group with splitting on one strand.''

The results in the paper are geometric and algebraic descriptions of
\(\widehat{BV}\) and \(BV\) that reveal their algebraic structure, a
derivation of a normal form for the elements of the groups, an
infinite presentation for \(\widehat{BV}\), and sketches of
arguments that the geometric and algebraic descriptions of each
group are of isomorphic groups.  All that we say applies with
trivial modification to \(V\), and we get a similar normal form for
\(V\).  This normal form for \(V\) is general knowledge but has
never been recorded.  In \cite{brin:bv3}, we derive finite
presentations for \(\widehat{BV}\) and \(BV\), and also a new finite
presentation for \(V\) that is similar to that of \(BV\).

Patrick Dehornoy \cite{dehornoy:geom-presentations} has
independently discovered \(BV\) and \(\widehat{BV}\).  Dehornoy's
motivation is not that removed from ours---he builds the groups as
braided versions of structure groups of algebraic identities---but
his techniques of construction, analysis, and set of observations
about the groups are different.  The group \(BV\) is the main focus
in the current paper, and we construct \(\widehat{BV}\) primarily to
get to \(BV\).  The focus in \cite{dehornoy:geom-presentations} is
on the group \(\widehat{BV}\) (called \(FB_\infty\) in that paper)
and its strong relationship with the law of left self
distributivity.  See \cite{MR2001j:20057} for more information and
consequences of the relationship between this law and the ordinary
braid groups.

The groups \(BV\) and \(\widehat{BV}\) are related to other groups
in the literature.  The group \(BV\) injects into the ``universal
mapping class group'' \(\scr{B}\) of \cite{math.GT/0210007}.  There
is also a non-finitely generated braided version of \(V\)
constructed in \cite{green+serg} that is different from \(BV\) and
\(\widehat{BV}\) since \(BV\) and \(\widehat{BV}\) have finite
presentations.  The group \(\widehat{V}\) could have been defined in
a single word the proof of Proposition 4.1 of
\cite{brown:finiteprop} (where it would have been \(G_{2,\infty}\)
in the notation of that paper).

\subsection{Descriptions of \protect\(V\protect\),
\protect\(BV\protect\) and
\protect\(\widehat{BV}\protect\)}\mylabel{GeomDescrSec} Elements of
\(V\) are most easily described using the standard ``deleted middle
thirds'' description of the Cantor set \(C\).  The set \(C\) is a
limit of a sequence of collections of closed intervals in the unit
interval \([0,1]\).  The first few collections are 
\mymargin{CoveringCollns}\begin{equation}\label{CoveringCollns}
\begin{split}
\{[0&,1]\} \\ \{[0,\hbox{\(\frac{1}{3}\)}]\,\,&,\,\,
[\hbox{\(\frac{2}{3}\)},1]\} \\
\{[0,\hbox{\(\frac{1}{9}\)}]\,\,,\,\, [\hbox{\(\frac{2}{9},
\frac{1}{3}\)}]\,\,&,\,\, [\hbox{\(\frac{2}{3},
\frac{7}{9}\)}]\,\,,\,\, [\hbox{\(\frac{8}{9}\)}, 1]\} \\ \{
\hbox{\( [0,\frac{1}{27}]\,\,,\,\, [\frac{2}{27},
\frac{1}{9}]\,\,,\,\, [\frac{2}{9}, \frac{7}{27}]\,\,,\,\,
[\frac{8}{27}, \frac{1}{3}]\)}\,\,&,\,\, \hbox{\(
[\frac{2}{3},\frac{19}{27}]\,\,,\,\,
[\frac{20}{27},\frac{7}{9}]\,\,,\,\, [\frac{8}{9},
\frac{25}{27}]\,\,,\,\, [\frac{26}{27}, 1]\)}\} \\ &\vdots
\end{split}\end{equation}

Elements of \(V\) are defined using covers of \(C\) by pairwise
disjoint intervals chosen from the collections in
\tref{CoveringCollns}.  Given a pair 
of such covers with the same number of intervals (necessarily
finite) and a one-to-one correspondence between them, an element of
\(V\) is created by taking each interval in the first cover affinely
in an orientation preserving way onto the corresponding interval in
the second cover and restricting this map to \(C\).  The restriction
is demonstrably a homeomorphism of \(C\).  The group \(V\) is the
set of all such maps.  Below, we show one such map \(f\) where we
indicate the bijection by numbering the intervals.  \[ \left(\,\,{\xy
(0,0); (13.5,0)**@{-}; (6.75,2)*{1}; (27,0); (32.5,0)**@{-};
(29.75,2)*{2}; (36,0); (40.5,0)**@{-}; (38.25,2)*{3}; 
\endxy}\,\,\right) \qquad {\xy (6,2)*{f}; (0,0);
(12,0)**@{-};(12,0)*\dir{>} \endxy} \qquad \left(\,\,{\xy 
(0,0); (4.5,0)**@{-}; (2.25,2)*{3}; (9,0); (13.5,0)**@{-};
(11.25,2)*{1}; (27,0); (40.5,0)**@{-}; (33.75,2)*{2}; \endxy}\,\,\right)
\] The map \(f\) is the restriction of the following affine
surjections \[ \hbox{\([0,\frac{1}{3}]\)} \rightarrow
\hbox{\([\frac{2}{9},\frac{1}{3}]\)}, \qquad \hbox{\([\frac{2}{3},
\frac{7}{9}]\)} \rightarrow \hbox{\([\frac{2}{3}, 1]\)}, \qquad
\hbox{\([\frac{8}{9},1]\)} \rightarrow \hbox{\([0, \frac{1}{9}]\)}
\] to the portions of \(C\) contained in the given intervals.

Elements of \(V\) are usually coded by pairs of labeled binary
trees.  For example, the map \(f\) above is coded by the pair \[
\left(
\xy
(0,0); (5,5)**@{-}; (10,0)**@{-}; (5,-5)**@{-};
(10,0); (15,-5)**@{-};
(0,-2)*{1}; (5,-7)*{2}; (15,-7)*{3};
\endxy
\quad, \quad
\xy
(0,0); (-5,5)**@{-}; (-10,0)**@{-}; (-5,-5)**@{-};
(-10,0); (-15,-5)**@{-};
(-15,-7)*{3}; (-5,-7)*{1}; (0,-2)*{2};
\endxy
\right).
\]
The structure of the left tree indicates that the interval \([0,1]\)
(corresponding to the root at the top of the tree) is to be split,
and the resulting right interval \([\frac{2}{3},1]\) is to be split
again.  This describes the intervals in the domain of \(f\).  The
right tree codes the splittings needed to describe
the intervals in the range of \(f\).

To obtain an element of \(BV\), we embed \(C\) in the plane \(\R^2\)
as a subset of the \(x\)-axis.  Let \(C\) be covered by a collection
\(A\) of pairwise disjoint intervals from \tref{CoveringCollns} and
also a similar collection \(B\) with the same number of intervals.
An element of \(BV\) will take intervals in \(A\) to intervals in
\(B\) exactly as described above, but the move will be accomplished
by an isotopy of \(\R^2\) with compact support.  That is, the move
will be accomplished by braiding if we  view the isotopy as a level
preserving homeomorphism from \(\R^2\times [0,1]\) to itself and
letting the braid strands be the images of the components of
\(C\times[0,1]\).  A restriction that must be observed is that
during the isotopy, each interval in \(A\) must have its image
during the isotopy parallel to the \(x\)-axis at all times.
Isotopies \(u\) and \(v\) are equivalent if there is a level
preserving isotopy of \(\R^2\times [0,1]\) from \(u\) to \(v\)
(adhering to the restriction that the images of intervals from \(A\)
be kept parallel to the \(x\)-axis throughout)  that are fixed on
the Cantor set at the 0 and 1 levels.  Thus \(BV\) is seen to be a
subgroup of a braid group on a Cantor set of strands.

The surjection from \(BV\) to \(V\) is obtained by taking each
element of \(BV\) to the homeomorphism of \(C\) obtained at the end
of the isotopy.

An element of \(BV\) can also be coded by pairs of binary trees, but
now the connection from the leaves of the first to the leaves of the
second is given by a braid and not a bijection.  This is most easily
pictured by drawing the second tree upside down below the first and
drawing the braid connecting the leaves between the two trees.

As an example, the following is one element of \(BV\) (out of
infinitely many) that maps to the element \(f\) of \(V\) in the
example above.  We draw both the ``trees and braid'' encoding of the
element as well as a picture of a braiding of intervals.
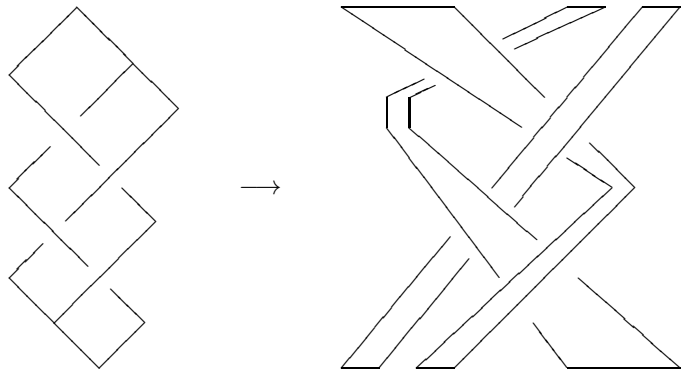
\[
\xy
(0,24); (-9,15)**@{-}; (-9,15)**@{-}; (-3,9)**@{-};
(3,3)**@{-}; (6,0); (10.5,-4.5)**@{-}; (-3,-18)**@{-};
(0,24); (7.5,16.5)**@{-}; 
(0.5, 9.5)**@{-}; (-3.5,5.5);
(-9,0)**@{-};
(1.5,-10.5)**@{-};
(4.5,-13.5); (9,-18)**@{-}; (3,-24)**@{-};
(7.5,16.5); (13.5,10.5)**@{-}; (-1.5,-4.5)**@{-}; (-4.5,-7.5);
(-9,-12)**@{-};
(3,-24)**@{-};
\endxy
\qquad\longrightarrow\qquad
\xy
(0,24); (15,24)**@{-};
(30,24); (35,24)**@{-};
(40,24); (45,24)**@{-};
(30,24); (21.5,19.75)**@{-};
(11, 14.5); (6,12)**@{-}; 
(6,8)**@{-}; 
(21, -12)**@{-};
(25.5, -18); 
(30,-24)**@{-};
(35,24); (23, 18.5)**@{-};
(12.5, 13.6); (9,12)**@{-}; 
(9,8)**@{-}; (26,-7)**@{-};
(31.5,-12);
(45,-24)**@{-};
(0,24); 
(24,8)**@{-};
(30,4); (36,0)**@{-}; (10,-24)**@{-};
(15,24); 
(27,12)**@{-};
(33,6); (39,0)**@{-}; (15,-24)**@{-};
(40,24); (20,0)**@{-}; (14.5,-6.5);  (0,-24)**@{-};
(45,24); (23,-2.5)**@{-}; (17,-9.5); (5,-24)**@{-};
(0,-24); (5,-24)**@{-};
(10,-24); (15,-24)**@{-};
(30,-24); (45,-24)**@{-};
\endxy
\]

The group \(\widehat{BV}\) is built in the same was as \(BV\) but it
acts on a countable number of copies of \(C\).  Specifically, we
embed one copy of \(C\) in the interval \([2i,\,2i+1]\) in the
\(x\)-axis for each integer \(i\ge0\) to create a subset \(X\) of
\(\R^2\).  Now infinite covers of \(X\) by pairwise disjoint closed
intervals are used where all but finitely many of the intervals are
intervals of the form \([2i,\,2i+1]\).  The remaining intervals in
each cover are translates of the intervals in \tref{CoveringCollns}
by even integers.  Again, an isotopy takes one cover to the other
keeping the images of the chosen intervals parallel to the
\(x\)-axis throughout.  Since the number of intervals is infinite,
there is no restriction that the number of intervals be equal.
However, this might result in a ``shift'' taking place for large
values of \(x\).  It is required that this shift be done by an
isotopy that, outside a compact set, consists of \((x,y)\mapsto
(x+td(1-|y|), y)\) for \(|y|<1\) and \(x\) greater than some \(K\)
and is the identity otherwise.  The integer \(d\) is the total
amount of shift and \(t\) is the parameter of the isotopy.
Isotopies between isotopies that create the equivalence classes of
the braids are required to have compact support, and, as above, are
required to keep the images of the defining intervals parallel to
the \(x\)-axis throughout.

We are in a position to explain some remarks made above.  The
structure of the trees is what keeps track of the ``order of
splitting.''  The restriction that the isotopies used by elements of
\(BV\) keep certain intervals parallel to the \(x\)-axis explains
why \(BV\) is a proper of the group \(\scr{B}\) of
\cite{math.GT/0210007} since the group \(\scr{B}\) allows the intervals
to rotate.

\subsection{Multiplication in \protect\(BV\protect\)} Because
strands are not really split and combined (bundles of strands over
the Cantor set are simply degrouped and regrouped), the following
relations are clear.  \mymargin{CancMoves} \begin{equation}
\label{CancMoves} \xy (-2,4);(0,2)**@{-};(0,-2)**@{-};
(-2,-4)**@{-}; (2,4);(0,2)**@{-};(2,-4);(0,-2)**@{-} \endxy
\leftrightarrow \xy (-2,4);(-2,-4)**@{-}; (2,4);(2,-4)**@{-} \endxy
\qquad\qquad\qquad\qquad \xy
(0,6);(0,4)**@{-};(-2,2)**@{-};(-2,-2)**@{-};(0,-4)**@{-};(0,-6)**@{-};
(0,4);(2,2)**@{-};(2,-2)**@{-};(0,-4)**@{-} \endxy \leftrightarrow
\xy (0,6);(0,-6)**@{-} \endxy \end{equation} The pictures in
\tref{CancMoves} simply  express
the fact that certain splitting and recombining operations are
inverses of each other.

We are in a position to multiply pictorally.  As in the braid group,
the product \(uv\) of \(u\) and \(v\) is drawn by putting \(u\) over
\(v\).  If \(u\) is the element 
\(
\xy
(0,6); (0,4)**@{-}; (2,2)**@{-}; (0,0)**@{-}; (-2,-2)**@{-}; 
(0,-4)**@{-}; (0,-6)**@{-};
(0,4); (-4,0)**@{-}; (-2,-2)**@{-};
(2,2); (4,0)**@{-}; (0,-4)**@{-};
\endxy
\), then \(u^2\) is calculated as follows.

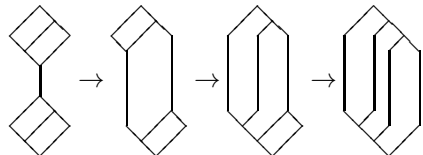
\[
\xy
(0,10);(-4,6)**@{-};(0,2)**@{-};(0,-2)**@{-};
(-4,-6)**@{-};(0,-10)**@{-};
(0,10);(4,6)**@{-};(0,2)**@{-};(0,-2)**@{-};
(4,-6)**@{-};(0,-10)**@{-};
(2,8);(-2,4)**@{-};  (2,-4);(-2,-8)**@{-}
\endxy
\rightarrow
\xy
(-1,10);(-5,6)**@{-};(-3,4)**@{-};
(-3,-6)**@{-};(1,-10)**@{-};
(-1,10);(3,6)**@{-};(3,-4)**@{-};
(5,-6)**@{-};(1,-10)**@{-};
(1,8);(-3,4)**@{-};  (3,-4);(-1,-8)**@{-}
\endxy
\rightarrow
\xy
(-1,10);(-5,6)**@{-};(-5,-4)**@{-};
(1,-10)**@{-};
(-1,10);(3,6)**@{-};(3,-4)**@{-};
(5,-6)**@{-};(1,-10)**@{-};
(1,8);(-1,6)**@{-};(-1,-4)**@{-};(-3,-6)**@{-};
(3,-4);(-1,-8)**@{-}
\endxy
\rightarrow
\xy
(-1,10);(-5,6)**@{-};(-5,-4)**@{-};
(1,-10)**@{-};
(-1,10);(5,4)**@{-};
(5,-6)**@{-};(1,-10)**@{-};
(1,8);(-1,6)**@{-};(-1,-4)**@{-};(-3,-6)**@{-};
(3,6);(1,4)**@{-};(1,-6)**@{-};(-1,-8)**@{-}
\endxy
\]

The above example has a trivial braid between two trees.
A product of three elements, two of which involve non-trivial braids
is shown below.
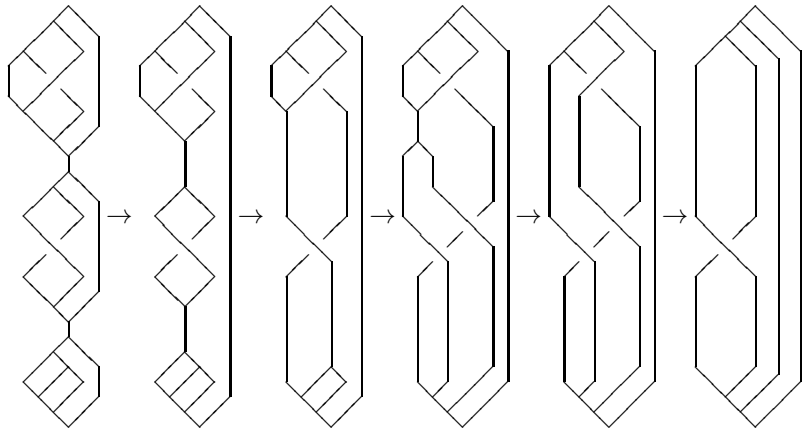
\[
\xy
(0,28);
(-8,20)**@{-};(-8,16)**@{-};
(0,8)**@{-};(0,6)**@{-};
(-6,0)**@{-};(2,-8)**@{-};(-2,-12)**@{-};(0,-14)**@{-};(0,-16)**@{-};
(-6,-22)**@{-};(0,-28)**@{-};
(0,28);(4,24)**@{-};(4,12)**@{-};(0,8)**@{-};
(0,6);(4,2)**@{-};(4,-10)**@{-};(0,-14)**@{-};
(0,-16);(4,-20)**@{-};(4,-24)**@{-};(0,-28)**@{-};
(-2,26);(2,22)**@{-};(-6,14)**@{-};
(-6,22);(-3,19)**@{-};(-1,17);(2,14)**@{-};(-2,10)**@{-};
(-2,4);(2,0)**@{-};(-1,-3)**@{-};(-3,-5);(-6,-8)**@{-};(-2,-12)**@{-};
(-2,-18);(2,-22)**@{-};(-2,-26)**@{-};
(0,-20);(-4,-24)**@{-};
\endxy
\rightarrow
\xy
(0,28);
(-8,20)**@{-};(-8,16)**@{-};
(-2,10)**@{-};(-2,4)**@{-};
(-6,0)**@{-};(2,-8)**@{-};(-2,-12)**@{-};(-2,-12)**@{-};(-2,-18)**@{-};
(-6,-22)**@{-};(0,-28)**@{-};
(0,28);(4,24)**@{-};(4,-24)**@{-};(0,-28)**@{-};
(-2,26);(2,22)**@{-};(-6,14)**@{-};
(-6,22);(-3,19)**@{-};(-1,17);(2,14)**@{-};(-2,10)**@{-};
(-2,4);(2,0)**@{-};(-1,-3)**@{-};(-3,-5);(-6,-8)**@{-};(-2,-12)**@{-};
(-2,-18);(2,-22)**@{-};(-2,-26)**@{-};
(0,-20);(-4,-24)**@{-};
\endxy
\rightarrow
\xy
(0,28);
(-8,20)**@{-};(-8,16)**@{-};
(-6,14)**@{-};(-6,0)**@{-};
(-6,0)**@{-};(0,-6)**@{-};(0,-20)**@{-};
(2,-22)**@{-};(-2,-26)**@{-};
(0,28);(4,24)**@{-};(4,-24)**@{-};(0,-28)**@{-};
(-2,26);(2,22)**@{-};(-6,14)**@{-};
(-6,22);(-3,19)**@{-};(-1,17);(2,14)**@{-};
(2,0)**@{-};(-1,-3)**@{-};(-3,-5);(-6,-8)**@{-};
(-6,-22)**@{-};(0,-28)**@{-};
(0,-20);(-4,-24)**@{-};
\endxy
\rightarrow
\xy
(0,28);
(-8,20)**@{-};(-8,16)**@{-};
(-6,14)**@{-};
(-6,10)**@{-};(-8,8)**@{-};(-8,0)**@{-};(-2,-6)**@{-};(-2,-22)**@{-};
(-4,-24)**@{-};
(-2,-26);(4,-20)**@{-};(4,-4)**@{-};(-4,4)**@{-};(-4,8)**@{-};(-6,10)**@{-};
(0,28);(6,22)**@{-};(6,-22)**@{-};(0,-28)**@{-};
(-2,26);(2,22)**@{-};(-6,14)**@{-};
(-6,22);(-3,19)**@{-};(-1,17);(4,12)**@{-};
(4,2)**@{-};(2,0)**@{-};(0,-2);(-2,-4)**@{-};(-4,-6);(-6,-8)**@{-};
(-6,-22)**@{-};(0,-28)**@{-};
\endxy
\rightarrow
\xy
(0,28);
(-8,20)**@{-};(-8,0)**@{-};(-2,-6)**@{-};(-2,-22)**@{-};
(-4,-24)**@{-};
(-2,-26);(4,-20)**@{-};(4,-4)**@{-};(-4,4)**@{-};
(-4,16)**@{-};(2,22)**@{-};(-2,26)**@{-};
(0,28);(6,22)**@{-};(6,-22)**@{-};(0,-28)**@{-};
(-6,22);(-3,19)**@{-};(-1,17);(4,12)**@{-};
(4,2)**@{-};(2,0)**@{-};(0,-2);(-2,-4)**@{-};(-4,-6);(-6,-8)**@{-};
(-6,-22)**@{-};(0,-28)**@{-};
\endxy
\rightarrow
\xy
(0,28);
(-8,20)**@{-};(-8,0)**@{-};(0,-8)**@{-};(0,-20)**@{-};
(-4,-24)**@{-};
(-2,-26);(3,-21)**@{-};(3,21)**@{-};(-2,26)**@{-};
(0,28);(6,22)**@{-};(6,-22)**@{-};(0,-28)**@{-};
(-4,24);(0,20)**@{-};
(0,0)**@{-};(-3,-3)**@{-};(-5,-5);(-8,-8)**@{-};
(-8,-20)**@{-};(0,-28)**@{-};
\endxy
\]

The restrictions that we impose make a difference.  Imagine
that the following moves are allowed in a part of a diagram.

\mymargin{ForbiddenMoves}
\begin{equation}
\label{ForbiddenMoves}
\xy
(0,6);(0,4)**@{-};(-4,0)**@{-};(0,-4)**@{-};(0,-6)**@{-};
(0,4);(4,0)**@{-};(0,-4)**@{-};
(-2,2);(2,-2)**@{-}
\endxy
\leftrightarrow
\xy
(0,6);(0,4)**@{-};(-4,0)**@{-};(0,-4)**@{-};(0,-6)**@{-};
(0,4);(4,0)**@{-};(0,-4)**@{-};
(2,2);(-2,-2)**@{-}
\endxy
\qquad\qquad\qquad
\xy
(3,3);(3,1)**@{-};(1,-1)**@{-};(1,-3)**@{-};
(-3,-3);(-3,-1)**@{-};(-1,1)**@{-};(-1,3)**@{-};
(1,-1);(-1,1)**@{-}
\endxy
\leftrightarrow
\xy
(-3,3);(-3,1)**@{-};(-1,-1)**@{-};(-1,-3)**@{-};
(3,-3);(3,-1)**@{-};(1,1)**@{-};(1,3)**@{-};
(-1,-1);(1,1)**@{-}
\endxy
\end{equation}

The move on the left can be accomplished by a 180 degree rotation inside
a 2-sphere that intersects the diagram in exactly two points.  This
move would correspond to an illegal rotation of an interval during
an isotopy.  The move on the right is accomplished by an isotopy
that reverses the slope of the line joining the two strands but that
does not interchange the two strands.  This move essentially alters
a tree that is part of the diagram.  Note that the first move can
also be realized as an application of the second move to the line in
the middle of the square.

The moves in \tref{ForbiddenMoves} have the following consequence.
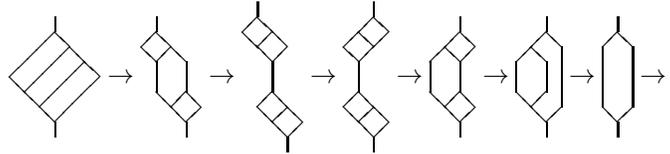
\[\xy
(0,8); (0,6)**@{-}; (-6,0)**@{-}; (-2,-4)**@{-}; (0,-6)**@{-};
(0,6); (2,4)**@{-}; (6,0)**@{-}; (0,-6)**@{-}; (0,-8)**@{-};
(2,4); (-4, -2)**@{-};
(4,2); (-2,-4)**@{-} 
\endxy
\rightarrow
\xy
(-2,8); (-2,6)**@{-}; 
(-4,4)**@{-};(-2,2)**@{-};(-2,-2)**@{-};(2,-6)**@{-};(2,-8)**@{-};
(2,-6);(4,-4)**@{-};(2,-2)**@{-};(2,2)**@{-};(0,4)**@{-};(-2,6)**@{-};
(-2,2);(0,4)**@{-};
(2,-2);(0,-4)**@{-}
\endxy
\rightarrow
\xy
(-2,10);(-2,8)**@{-};(-4,6)**@{-};(0,2)**@{-};(0,-2)**@{-};
(4,-6)**@{-};(2,-8)**@{-};(2,-10)**@{-};
(-2,8);(2,4)**@{-};(0,2)**@{-}; (0,6);(-2,4)**@{-};
(2,-8);(-2,-4)**@{-};(0,-2)**@{-}; (0,-6);(2,-4)**@{-}
\endxy
\rightarrow
\xy
(2,10);(2,8)**@{-};(4,6)**@{-};(0,2)**@{-};(0,-2)**@{-};
(4,-6)**@{-};(2,-8)**@{-};(2,-10)**@{-};
(2,8);(-2,4)**@{-};(0,2)**@{-}; (0,6);(2,4)**@{-};
(2,-8);(-2,-4)**@{-};(0,-2)**@{-}; (0,-6);(2,-4)**@{-}
\endxy
\rightarrow
\xy
(2,8);(2,6)**@{-};(4,4)**@{-};(2,2)**@{-};(2,-2)**@{-};
(4,-4)**@{-};(2,-6)**@{-};(2,-8)**@{-};
(2,6);(-2,2)**@{-};(-2,-2)**@{-};(2,-6)**@{-};
(0,-4);(2,-2)**@{-}; (0,4);(2,2)**@{-}
\endxy
\rightarrow
\xy
(2,8);(2,6)**@{-};(4,4)**@{-};
(4,-4)**@{-};(2,-6)**@{-};(2,-8)**@{-};
(2,6);(-2,2)**@{-};(-2,-2)**@{-};(2,-6)**@{-};
(0,-4);(2,-2)**@{-};(2,2)**@{-};(0,4)**@{-}
\endxy
\rightarrow
\xy
(2,8);(2,6)**@{-};(4,4)**@{-};
(4,-4)**@{-};(2,-6)**@{-};(2,-8)**@{-};
(2,6);(0,4)**@{-};(0,-4)**@{-}; (2,-6)**@{-};
\endxy
\rightarrow
\xy
(2,8);(2,-8)**@{-};
\endxy
\]
However, the figure on the left represents a non-trivial element of
\(BV\) and of \(V\).

The goal of the paper is to expose as much of the algebraic
structure of \(BV\) as possible.  For this reason, our construction
and analysis of \(BV\) is more algebraic than geometric.  The next
section reviews some of the techniques that we will use in the
remainder of the paper and contains a bit of an outline as to how
the techniques will be put to use.

\section{Conventions, definitions and needed
constructions}\mylabel{backgroundSec}

\subsection{Outline} Working backwards from the group \(BV\), we get
\(BV\) as a subgroup of the better behaved group \(\widehat{BV}\).
The group \(\widehat{BV}\) is obtained as a group of right fractions
of a cancellative monoid \(\FB\) and a normal form (reduced terms)
for the elements of \(\widehat{BV}\) is obtained from the fact that
pairs of elements of \(\FB\) have unique greatest common right
factors.  The monoid \(\FB\) is a \ZS product (a generalization of
the semidirect product in which neither factor need be normal) of
the monoid \(\scr{F}\) and the group \(B_\infty\).  The group
\(B_\infty\) is the familiar braid group on infinitely many strands.
The monoid \(\scr{F}\) is a monoid of binary forests and is
understood by obtaining a normal form for its elements.  The normal
form is obtained by a standard technique using the concept of
terminating and confluent relations---a technique that is used more
than once in this paper.

In the rest of this section, we first review the techniques for
getting normal forms from confluent and terminating relations.  Next
we review groups of fractions and the use of greatest common right
factors to obtain reduced terms.  Lastly, we review \ZS products.
Since we will need information about these products that lead to
certain properties (such as the existence of greatest common
factors) and also lead to presentations, we review what is
needed to get such information from the products.

\subsection{Distinguished representatives}

Normal forms will be used to establish properties of groups with
presentations.  The normal forms will come from standard techniques
from string rewriting which arise in turn from properties of
relations.  

A binary relation \(\rightarrow\) on a set \(A\) is called {\itshape
terminating} if there is no infinite sequence \(x_0\rightarrow
x_1\rightarrow \cdots\).  Note that a reflexive relation cannot be
terminating.  An element \(a\in A\) is said to be {\itshape
irreducible} in \(A\) if \(a\rightarrow x\) is false for all \(x\in
A\).

We let \(\overset{*}\rightarrow\) denote the reflexive, transitive
closure of \(\rightarrow\).  The relation \(\rightarrow\) is
{\itshape locally confluent} if for every \(x\), \(y\) and \(z\)
satisfying \(x\rightarrow y\) and \(x\rightarrow z\), there is a
\(w\) satisfying \(y \overset{*}\rightarrow w\) and \(z
\overset{*}\rightarrow w\).  If we let \(\sim\) denote the
equivalence relation generated by \(\rightarrow\) (the symmetric,
reflexive, transitive closure of \(\rightarrow\)), then we get the
following result of Newman \cite{newman:comb}.

\begin{prop}\mylabel{NewmanLemma} If a binary relation
\(\rightarrow\) on a set \(A\) is terminating and locally confluent,
then every equivalence class under \(\sim\) contains a unique
element that is irreducible in \(A\).  Further, \(x
\overset{*}\rightarrow a\) for every \(x\) where \(a\) is the unique
irreducible in \(A\) that is in the equivalence class containing
\(x\).   \end{prop}

A terminating, binary relation that is also locally confluent is
called {\itshape complete}.

Binary relations will often come from rewriting rules.  If
\(\Sigma\) is a set (which we call an alphabet in this situation),
then \(\Sigma^*\) denotes the set of strings (finite sequences) of
elements of \(A\).  The empty string (sequence of length zero) is
also in \(\Sigma^*\).  This is a monoid under concatenation and we
will refer to it as the {\itshape free monoid} on \(\Sigma\).  A
binary relation \(\rightarrow\) on \(\Sigma^*\) can be referred to
as a {\itshape rewriting rule} which terminology implies that
another binary relation \(\theta\) on \(\Sigma^*\) is to be regarded
as a consequence of \(\rightarrow\) as follows.  If \(u\), \(v\),
\(p\) and \(q\) are in \(\Sigma^*\), then we write \(puq \,\theta \,
pvq\) if \(u\rightarrow v\).  Confusingly, \(\rightarrow\) is often
used for both the rewriting rule and its consequence.  The confusion
is usually not crippling.

A rewriting rule is called {\itshape complete} if its consequence is
complete.  

If \(\langle X\mid R\rangle\) is a presentation, then we regard
\(R\) as a relation on words in \(X\).  As a relation, \(R\) is
usually thought of as either symmetric or its symmetric closure is
implied.  However, if the symmetric closure is not taken, then \(R\)
can also be thought of as a set of rewriting rules.  If it turns out
that \(R\) is complete as a set of rewriting rules, then we say that
the presentation is {\itshape complete}.  The power of a complete
presentation is that Proposition \ref{NewmanLemma} gives each
element of the presented object a distinguished representative.

\subsection{Properties of monoids and the Ore theorem }

The group of fractions construction will be important.  We set out
the necessary definitions and results.

A {\itshape semigroup} is a set with a binary, associative product.
A {\itshape monoid} is a semigroup with a global, two-sided
identity.  A {\itshape group} is a monoid in which every element has
a two-sided inverse.

A semigroup is {\itshape left cancellative} if \(ax=ay\) always
implies \(x=y\) and is {\itshape right cancellative} if \(xa=ya\)
always implies \(x=y\).  A semigroup is {\itshape cancellative} if
it is both right and left cancellative.  A semigroup is {\itshape
strongly left cancellative} if it is left cancellative and if
\(ab=a\) implies that \(b\) is a global, two-sided identity.  The
reader can define strongly right cancellative.

\begin{lemma}\mylabel{GetStrLftCanc} Let \(S\) be a cancellative
semigroup.  Then (1) \(S\) is strongly left cancellative and
strongly right cancellative and (2) if \(S\) has a global, two-sided
identity 1, and \(ab=1\), then \(ba=1\).  \end{lemma}

In a semigroup \(S\), a {\itshape right multiple} of an element
\(x\in S\) is an element \(y\in S\) for which there is an element
\(p\in S\) so that \(y=xp\).  A subset \(C\) of \(S\) has a
{\itshape common right multiple} \(z\) if \(z\) is a right multiple
of every element of \(C\).  Usually, \(C\) has two elements and we
refer to the common right multiple of \(C\) as the common right
multiple of the two elements.  A common right multiple \(z\) for
\(C\) is a {\itshape least common right multiple} if every common
right multiple for \(C\) is a right multiple of \(z\).  A semigroup
{\itshape has common right multiples} if every pair of elements has
a common right multiple.  A semigroup {\itshape has least common
right multiples} if every pair of elements with a common right
multiple has a least common right multiple.  Note that the last
definition has been carefully worded so that a semigroup with least
common right multiples need not have common right multiples.

In the previous paragraph, every appearance of the word right can be
replaced by the word left to give a corresponding set of
definitions.  It will turn out that the important concepts for this
paper will be ``common right multiples'' and ``least common left
multiples.''  The first concept will lead to groups of fractions and
the second concept will lead to distinguished representatives in the
group of fractions.

Least common multiples are often associated with greatest common
factors.  An element \(r\) in a semigroup \(S\) is a {\itshape right
factor} of an element \(x\in S\) if there is a \(p\in S\) so that
\(x=pr\).  Two elements \(x\) and \(y\) in \(S\) have a {\itshape
common right factor} \(r \in S\) if \(r\) is a right factor of both
\(x\) and \(y\).  The common right factor \(r\) is a {\itshape
greatest common right factor} if every common right factor of \(x\)
and \(y\) is a right factor of \(r\).  A semigroup \(S\) {\itshape
has greatest common right factors} if every pair of elements with a
common right factor has a greatest common right factor.

A {\itshape length function} on a semigroup \(S\) is a homomorphism
to the natural numbers \(\N\) so that the preimage of 0 is contained
in the set of those \(x\in S\) for which there is a \(y\in S\) so
that \(xy=yx\) is a global, two-sided identity for \(S\).

The following is a pleasant exercise for the reader.

\begin{lemma}\mylabel{MaxEquivs} A semigroup with a length function
has least common left multiples if and only if it has greatest
common right factors.  \end{lemma}

In the following, a presentation \(\langle X\mid R\rangle\) is
thought of as a set \(X\) and a relation \(R\) on either the free
monoid \(X^*\) on \(X\) if the presentation is a monoid presentation
or the free group \(F_X\) on \(X\) if the presentation is a group
presentation.  As a relation, \(R\) is thought of as a set of
ordered pairs.  The two entries an ordered pair in \(R\) are
regarded as equal in the object being presented.

We now discuss groups of fractions.  Let \(S\) be a cancellative
semigroup with common right multiples.  Let pairs in \(S\times S\)
be denoted by \(\frac{x}{y}\) instead of \((x,y)\).  Let \(\sim\) be
the equivalence relation generated by \(\frac{xz}{yz}\rightarrow
\frac{x}{y}\) and define a product on representatives by
\(\frac{x}{y}\frac{y}{z} = \frac{x}{z}\).  It follows from the
existence of common right multiples that any two classes have
representatives that can be multiplied.  The product is well defined
and we have the following.

\begin{prop}[Ore]\mylabel{Ore} Let \(S\) be a cancellative semigroup
with common right multiples.  Then the following hold.
{\claimenum \item The multiplication above turns the
equivalence classes on \(S\times S\) into a group \(G\).  \item If
\(a\) is a fixed element of \(S\), then sending \(x\in S\) to
\(\frac{ax}{z}\) is a homomorphic embedding of \(S\) into \(G\).
\item Every element of \(G\) is representable in the form
\(pn^{-1}\) where \(p\) and \(n\) are in the image of the
embedding in (b).  \item If \(\langle X\mid Y\rangle\) is a
semigroup presentation of \(S\), then \(\langle X\mid Y\rangle\) is
a group presentation of \(G\).  \end{enumerate}} \end{prop}

Items (a)--(c) of Proposition \ref{Ore} comprise a mirror image of
Theorem 1.23 of \cite{cliff+prest:I} where the existence of common
left multiples (there called right reversible) is used instead.
Item (d) is well known and straightforward.

If \(U\) is a cancellative semigroup with common right multiples,
and if \(i:U\rightarrow G\) is an injective homomorphism into a
group so that every element of \(G\) is of the form \(uv^{-1}\) with
both \(u\) and \(v\) in the image of \(i\), then we call \(i\) an
{\itshape Ore embedding} of \(U\) into a {\itshape group of right
fractions} of \(U\) and write \(\frac{u}{v}\) for \(uv^{-1}\).
Theorem 1.25 of \cite{cliff+prest:I} justifies calling \(G\)
\emph{the} group of right fractions of \(U\).

We get distinguished representatives (corresponding to reduced
fractions) in the group of fractions under certain circumstances.
To state the next lemma, we make a definition.  Let \(i:U\rightarrow
G\) be an embedding of a semigroup into a group of right fractions.
We say that a representative \(\frac{u}{v}\) as above is {\itshape
in reduced terms} if whenever \(\frac{w}{z}\) is another such
representative with \(\frac{u}{v}\sim\frac{w}{z}\), then there is an
\(x\) in the image of \(i\) with \(w=ux\) and \(z=vx\).  It follows
that if \(\frac{u}{v}\) and \(\frac{u'}{v'}\) are two
representatives of the same element and both are in reduced terms,
then there are \(x\) and \(y\) in the image of \(i\) so that
\(u'=ux\), \(u=u'y\), \(v'=vx\) and \(v=v'y\).  This implies that
\(u=uxy\) and \(u'=u'yx\).  We identify \(U\) with \(i(U)\) and use
Lemma \ref{GetStrLftCanc} to conclude that \(xy\) and \(yx\) are
both global, two-sided identities for \(U\).  That is, both \(x\)
and \(y\) are invertible.  To give brief terminology to this
situation, we say that the representatives \(\frac{u}{v}\) and
\(\frac{u'}{v'}\) {\itshape differ by invertible elements} of \(U\).

\begin{lemma}\mylabel{ReducedTerms} Let \(U\) be a cancellative
semigroup with common right multiples, least common left multiples
and a length function.  Let \(i:U\rightarrow G\) be an Ore embedding
of \(U\) into a group of right fractions.  Then every element of
\(G\) has a representative in reduced terms and any two
representatives in reduced terms of one element of \(G\) differ by
invertible elements of \(U\).  In particular, any representative
\(uv^{-1}\) of an element \(g\in G\) with \(u\) and \(v\) in
\(i(U)\) and the length of \(u\) minimal among such representatives
of \(g\) is in reduced terms.  \end{lemma}

\begin{proof} This can be given a short direct proof, taking the
last sentence as a starting point. \end{proof}

\subsection{\protect\ZS products}

In spite of the fact that the previous discussion mentioned
semigroups repeatedly, we will work entirely with monoids.  \ZS
products work in even greater generality than semigroups (see
\cite{brin:zs}), but we will restrict our discussion to monoids for
practical reasons.

The \ZS product is a generalization of the semidirect product.  It
is a generalization in that no normality is required.  Proofs of the
statements in this section and a short history of the product can be
found in \cite{brin:zs}.  We first show how the ingredients of the
product arise.  We start with groups to give a familiar setting and
then generalize to monoids.

Let \(G\) be a group with identity \(1\), and with subgroups \(U\)
and \(A\) satsifying \(U\cap A=\{1\}\) and \(G=UA\).  Then each
\(g\in G\) is uniquely expressible as \(g=u\alpha\) with \(u\in U\)
and \(\alpha\in A\).  With \(u\in U\) and \(\alpha\in A\), consider
\(\alpha u\in G\).  There are uniqe elements \(u'\in U\) and
\(\alpha'\in A\) so that \(\alpha u=u'\alpha'\).  This defines two
functions \((\alpha,u)\mapsto \alpha^u\in A\) and
\((\alpha,u)\mapsto \alpha\cdot u\in U\) that are unique in that
they satisfy \(\alpha u=(\alpha\cdot u)(\alpha^u)\) for all \(u\in
U\) and \(\alpha\in A\).  We now move to monoids.

\begin{lemma}\mylabel{ExistMutualActs} Let \(M\) be a monoid and let
\(U\) and \(A\) be submonoids of \(M\).  Assume that every \(x\in
M\) is uniquely expressible in the form \(x=u\alpha\) with \(u\in
U\) and \(\alpha\in A\).  Then there are functions \((A\times
U)\rightarrow A\) written \((\alpha,u)\mapsto \alpha^u\), and
\((A\times U)\rightarrow U\) written \((\alpha,u)\mapsto \alpha\cdot
u\) defined by the property that \(\alpha u=(\alpha\cdot
u)(\alpha^u)\).  \end{lemma}

In the setting of the above lemma, the functions \((\alpha,
u)\mapsto \alpha^u\) and \((\alpha, u)\mapsto \alpha\cdot u\)
defined on \(A\times U\) will be called the {\itshape mutual
actions defined by the multiplication}.  These actions are
``internally'' generated by the multiplication.  We can also impose
actions ``externally.''

Let \(U\) and \(A\) be monoids.  Assume there are
functions \(A\times U \rightarrow A\) written \((\alpha ,u)\mapsto
\alpha \cdot u\) and \(A\times U\rightarrow U\) written \((\alpha
,u)\mapsto \alpha ^u\).  We will call such functions {\itshape
mutual actions} between \(A\) and \(U\).

We now define a multiplication on \(U\times A\)
by\mymargin{TheMult}\begin{equation}\label{TheMult}(u,\alpha
)(v,\beta ) = (u(\alpha \cdot v), \alpha ^v\beta ).\end{equation}
This multiplication is well defined.  The following ties it to the
setting of Lemma \ref{ExistMutualActs} and is essentially Lemma 3.9
of \cite{brin:zs}.

\begin{lemma}\mylabel{ZappaRecon} Let \(M\) be a monoid and let \(U\)
and \(A\) be submonoids of \(M\).  Assume that every \(x\in M\) is
uniquely expressible in the form \(x=u\alpha\) with \(u\in U\) and
\(\alpha\in A\), and let \((\alpha,u)\mapsto \alpha^u\) and
\((\alpha,u)\mapsto \alpha\cdot u\) defined on \(A\times U\) be the
mutual actions defined by the multiplication.  Use these mutual
actions and \tref{TheMult} to build a multiplication on \(U\times
A\).  Then sending \((u,\alpha)\) in \(U\times A\) to \(u\alpha\) in
\(M\) is an isomorphism of monoids.  \end{lemma}

Assuming the setting, hypotheses and notation of Lemma
\ref{ZappaRecon}, we say that \(M\) is the {\itshape (internal) \ZS
product} of \(U\) and \(A\) and write \(M=U\zapprod A\).

In the setting of monoids, the mutual actions are not arbitrary.
The next lemma gives sample properties that they satisfy.  The
statement that \(U\) is a submonoid of \(M\) contains the assumption
that the identity for \(U\) is the identity for \(M\).  The
Lemma combines Lemma 3.2 and Corollary 3.3.1 from
\cite{brin:zs}.

\begin{lemma}\mylabel{MutualActs} Let \(M\) be a monoid and let
\(U\) and \(A\) be submonoids of \(M\).  Assume that every \(x\in
M\) is uniquely expressible in the form \(x=u\alpha\) with \(u\in
U\) and \(\alpha\in A\), and let \((\alpha,u)\mapsto \alpha^u\) and
\((\alpha,u)\mapsto \alpha\cdot u\) defined on \(A\times U\) be the
mutual actions defined by the multiplication.  Let \(\alpha\) and
\(\beta\) come from \(A\) and \(u\) and \(v\) come from \(U\).  Let
\(1_U\) and \(1_A\) denote the identities of \(U\) and \(A\),
respectively.  Then the following hold. {\claimenum \item
\((\alpha\beta)\cdot u = \alpha\cdot (\beta\cdot u)\).  \item
\((\alpha\beta)^u = \alpha^{(\beta\cdot u)}\beta^u\).  \item
\(\alpha\cdot (uv) = (\alpha\cdot u)(\alpha^u\cdot v)\).  \item
\(\alpha^{(uv)} = (\alpha^u)^{v}\).  \item \(\alpha^{1_U}=\alpha\).
\item \(1_A\cdot u=u\).  \item \(\alpha\cdot 1_u=1_u\).  \item
\((1_A)^u=1_A\).  \claimenumend} \end{lemma}

The properties of Lemma \ref{MutualActs} are exactly those needed to
define an external \ZS product of monoids.  The next lemma is item
(xv) of Lemma 3.13 of \cite{brin:zs}.

\begin{lemma}\mylabel{PartMulPropsII} Let \(U\) and \(A\) be monoids
with mutual actions \((\alpha ,u)\mapsto \alpha \cdot u\) and
\((\alpha ,u)\mapsto \alpha ^u\) defined on \(A\times U\).  Assume
(a)--(h) of Lemma \ref{MutualActs}.  Then the multiplication
\tref{TheMult} makes \(U\times A\) a monoid \(M\).  Further
\(\alpha\mapsto (1_U, \alpha)\) and \(u\mapsto(u, 1_A)\) are
homomorphic embeddings of \(A\) and \(U\), respectively, into \(M\)
so that \(M\) is the internal \ZS product of the images.
\end{lemma}

We can refer to the monoid \(M\) of Lemma \ref{PartMulPropsII} as
the {\itshape (external) \ZS product} of \(U\) and \(A\) and again
write \(M=U\zapprod A\).

We need to know when \ZS products have certain properties.  In order
to discuss this, we need to look at mutual actions as families of
functions.  If a function \((\alpha, u)\mapsto \alpha^u\) is defined
from \(A\times U\) to \(A\), then we think of this as a family of
functions from \(A\) to itself parametrized by elements of \(U\).
We say that \((\alpha, u)\mapsto \alpha^u\) forms a {\itshape
surjective family} of functions if for every \(u\in U\) and
\(\alpha\in A\) there is a \(\beta\in A\) so that
\(\alpha=\beta^u\).  We say that \((\alpha, u)\mapsto \alpha^u\)
forms an {\itshape injective family} of functions if
\(\alpha^u=\beta^u\) always implies that \(\alpha=\beta\).  A family
is {\itshape coconfluent} if whenever \(\alpha^u=\beta^v\), there
are \(\gamma\), \(p\) and \(q\) so that \(\alpha=\gamma^p\) and
\(\beta=\gamma^q\).  A family satisfying (d) of Lemma
\ref{MutualActs} is {\itshape strongly coconfluent} if whenever
\(\alpha^u=\beta^v\) and \(u\) and \(v\) have a common left
multiple, there are \(\gamma\), \(p\) and \(q\) so that
\(\alpha=\gamma^p\), \(\beta=\gamma^q\) and \(pu=qv\).  Similary
definitions can be made for \((\alpha, u)\mapsto \alpha\cdot u\)
defined from \(A\times U\) to \(U\).  

The following is Lemma 3.6 of \cite{brin:zs}.

\begin{lemma}\mylabel{StrongCoConf} Let \(A\times B\rightarrow A\)
written \((\alpha,u)\mapsto \alpha^u\) be strongly coconfluent.
Assume that \(B\) is a right cancellative semigroup, assume that
\(\alpha^u=\beta^v\) and assume that \(u\) and \(v\) have a least
common left multiple \(l=au=bv\).  Then there is a \(\gamma\in A\)
so that \(\alpha=\gamma^a\) and \(\beta=\gamma^b\).  \end{lemma}

The following comprises items (viii) and (ix) of Lemma 3.12 of
\cite{brin:zs} where the proof is left to the reader.

\begin{lemma}\mylabel{PartMulProps} Assume the notation, hypotheses
and conclusion of Lemma \ref{PartMulPropsII}.  Then the following
hold.

\begin{enumerate}

\item\mylabel{CancItem} If \(U\) and \(A\) are both right
cancellative and \((\alpha, u)\mapsto \alpha^u\) is an injective
family, then \(U\zapprod A\) is right cancellative.

\item\mylabel{CRMItem} If \(U\) and \(A\) both have common right
multiples and \((\alpha, u)\mapsto \alpha \cdot u\) is a surjective
family, then \(U\zapprod A\) has
common right multiples.
\end{enumerate}
\end{lemma}

Least common left multiples are a bit more complicated.  The
following is Lemma 3.14 of \cite{brin:zs}.

\begin{lemma}\mylabel{PartMulPropsLCLM} Assume the hypotheses,
notation and conclusion of Lemma \ref{PartMulPropsII}.  If \(U\) is
cancellative with least common left multiples, if \(A\) is a group,
and if \((\alpha, u)\mapsto \alpha^u\) is strongly coconfluent, then
\(M=U\zapprod A\) has least common left multiples.  Further, the
least common left multiple
\((r,\alpha)(u,\theta)=(s,\beta)(v,\phi)\) of \((u,\theta)\) and
\((v,\phi)\) in \(U\zapprod A\) can be constructed so that
\(r(\alpha\cdot u)=s(\beta\cdot v)\) is the least common left
multiple of \((\alpha\cdot u)\) and \((\beta\cdot v)\) in \(U\).  If
\(M\) is cancellative (e.g., \((\alpha, u)\mapsto \alpha^u\) is an
injective family), then any least common left multiple
\((r,\alpha)(u,\theta)=(s,\beta)(v,\phi)\) of \((u,\theta)\) and
\((v,\phi)\) in \(U\zapprod A\) has the property that
\(r(\alpha\cdot u)=s(\beta\cdot v)\) is the least common left
multiple of \((\alpha\cdot u)\) and \((\beta\cdot v)\) in \(U\).
\end{lemma}

We copy from \cite{brin:zs} some very specialized results about
presentations of \ZS products that fit the needs of this paper.

Assume that presentations \(\langle X \mid R\rangle\) and \(\langle
Y \mid T\rangle\) of monoids \(U\) and \(A\), respectively, are
given with \(X\cap Y=\emptyset\), and that functions \(Y\times
X\rightarrow Y^*\) written \((\alpha,u)\mapsto \alpha^u\) and
\(Y\times X\rightarrow X\) written \((\alpha,u)\mapsto \alpha\cdot
u\) are given.  The unequal treatment of the codomains (\(X\) in one
case and \(Y^*\) in the other) is deliberate.

We extend these to functions \(Y^*\times X^*\rightarrow Y^*\) and
\(Y^*\times X^*\rightarrow X^*\) as follows.  Form the monoid
presentation
\mymargin{ZappaPresInt}\begin{equation}\label{ZappaPresInt}\langle
X\cup Y \mid Z\rangle\end{equation} in which \(Z\) is regarded as a
set of rewriting rules and consists of all pairs \((\alpha
u\rightarrow (\LeftAct{\alpha}{u})(\alpha^u))\) for \((\alpha,u)\in
Y\times X\).  The following easy lemma is Lemma 3.18 of
\cite{brin:zs}.

\begin{lemma}\mylabel{ZapPresIntCompl} The presentation
\tref{ZappaPresInt} is complete.  \end{lemma}

The irreducibles are of the form \(u\alpha\) with \(u\) a word in
the alphabet \(X\) and \(\alpha\) a word in the alphabet \(Y\).
This expresses the monoid presented by \tref{ZappaPresInt} as a \ZS
product of the free monoids \(X^*\) and \(Y^*\).  From Lemma
\ref{MutualActs} we get our desired extensions and the fact that
they satisfy the conclusions of that lemma.

The following combines Lemmas 3.17 and 3.19 of \cite{brin:zs}.

\begin{lemma}\mylabel{TwistedPresIII} Assume that presentations
\(\langle X \mid R\rangle\) and \(\langle Y \mid T\rangle\) of
monoids \(U\) and \(A\), respectively, are given with \(X\cap
Y=\emptyset\), and that functions \(Y\times X\rightarrow Y^*\)
written \((\alpha,u)\mapsto \alpha^u\) and \(Y\times X\rightarrow
X\) written \((\alpha,u)\mapsto \alpha\cdot u\) are given.  Let
\(\sim_R\) and \(\sim_T\) denote the equivalence relations on
\(X^*\) and \(Y^*\), respectively, imposed by the relation sets
\(R\) and \(T\), respectively.

Let the functions be extended to \(Y^*\times X^*\) as above and
assume that they satisfy the following.  If \((u,v)\) is in \(R\),
then for all \(\alpha\in Y\) we have \((\alpha\cdot u, \alpha\cdot
v)\) or \((\alpha\cdot v, \alpha\cdot u)\) is in \(R\) and
\(\alpha^u\sim_T\alpha^v\).  If \((\alpha,\beta)\) is in \(T\), then
for all \(u\in X\) we have \(\alpha\cdot u=\beta\cdot u\) and
\(\alpha^u\sim_T\beta^u\).  Then the extensions induce well defined
functions \(A\times U\rightarrow A\) and \(A\times U\rightarrow U\)
that satisfy the hypotheses (and thus the conclusions) of Lemma
\ref{PartMulPropsII} and the restriction of the function \(A\times
U\rightarrow U\) to \(A\times X\) has its image in \(X\).  Further a
presentation for the structure defined on \(U\zapprod A\) is
\mymargin{ZappaPresII}\begin{equation}\label{ZappaPresII}\langle
X\cup Y \mid R\cup T\cup W\rangle\end{equation} in which \(W\)
consists of all pairs \((\alpha u,(\LeftAct{\alpha}{u})(\alpha^u))\)
for \((\alpha,u)\in Y\times X\).  \end{lemma}

\section{The monoid of forests}

Our first algebraic structure will be a monoid whose objects are
forests.  For us a forest is a sequence of finite trees, only
finitely many of which are non-trivial.  We now give detailed
definitions.  None are surprising, but we give details to bring
reader and author into agreement on terminology.

The complete binary tree \(\scr{T}\) is the set of all finite
sequences with values in the set \(\{0,1\}\).  The sequence of
length 0, denoted \(\phi\), is included.  The sequences will be
referred to as strings, and we will concatenate string \(\alpha\)
and string \(\beta\) to give the string \(\alpha\beta\) in which
\(\alpha\) comes first and \(\beta\) comes last.  The most important
relation is prefix defined by ``\(\alpha\) is a prefix of
\(\alpha\beta\).''  The transitive closure of proper prefix is
ancestor and the inverse of ancestor is descendent.  The children of
\(u\) are exactly \(u0\) and \(u1\).

A finite binary tree \(T\) is a non-empty subset of \(\scr{T}\) that
is closed under ancestor and for which \(u0\) is in \(T\) if and
only if \(u1\) is in \(T\).  Every tree in this paper except
\(\scr{T}\) will be finite and binary, so we will stop using those
words as adjectives for trees.

Every tree \(T\) includes \(\phi\) which is called the root of
\(T\).  Elements of \(T\) will be called nodes, and the leaves of
\(T\) are the nodes of \(T\) with no children.  A tree is
non-trivial if it has more than one node.

We define trees this way so that if \(T\) and \(U\) are trees, then
\(T\cup U\) and \(T\cap U\) make sense.  It is elementary that both
\(T\cup U\) and \(T\cap U\) are trees.

There is a unique total ordering of the nodes of a tree \(T\) so
that every triad \(\{u, u0, u1\}\) in \(T\) is ordered \(u0<u<u1\).
We call this the left-right ordering of \(T\).  The restriction of
this order to the leaves of \(T\) is the left-right ordering of the
leaves of \(T\).

A forest \(F\) is an infinite sequence (function with domain \(\N\))
of trees so that all but finitely many are trivial.  We write
\(F_i\) for the \(i\)-th tree in \(F\).  The set \(\scr{F}\) of all
forests will be endowed shortly with a binary operation.

If \(v\) is a node of \(F_i\) for a forest \(F\), we distinguish it
from nodes of other trees in \(F\) by writing \(i.v\).  We order all
the leaves of \(F\) by giving the leaves of each tree in \(F\) the
left-right ordering and then insisting that all leaves in \(F_i\)
come before all the leaves in \(F_j\) when \(i<j\).

We number the leaves of \(F\) by the unique order preserving
function from the leaves to \(\N\).  The roots are numbered in the
obvious way: the \(i\)-th root is the root of \(F_i\).

If \(F\) and \(G\) are two forests, then we form \(FG\) by
identifying the \(i\)-th root of \(G\) with the \(i\)-th leaf of
\(F\).  Defining \(v_jG_j\) to mean \(\{v_ju\mid u\in G_j\}\) where
\(v_j\) is the \(j\)-th leaf of \(F\), then we can formally define
the \(i\)-th tree of \(FG\) as the union of \(F_i\) with all
\(v_jG_j\) where \(v_j\) is a leaf of \(F_i\).

Below we give
an example of a product \(FG\) of forests \(F\) and \(G\).  For
clarity, we have numbered the leaves of \(F\) and the roots of
\(G\) and \(FG\).

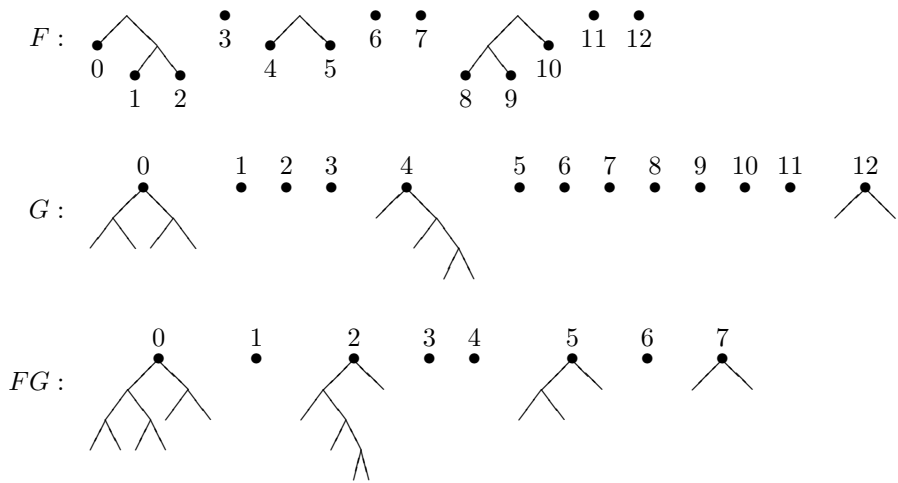
\[\begin{split}
F:\quad
&\xy
(0,0);(4,4)**@{-};(8,0)**@{-};(11,-4)**@{-};(8,0);(5,-4)**@{-};
(0,0)*{\bullet};(5,-4)*{\bullet};(11,-4)*{\bullet};
(0,-3)*\txt{0};(5,-7)*\txt{1};(11,-7)*\txt{2};
(17,4)*{\bullet};(17,1)*\txt{3};
(23,0);(27,4)**@{-};(31,0)**@{-};
(23,0)*{\bullet};(31,0)*{\bullet};
(23,-3)*\txt{4};(31,-3)*\txt{5};
(37,4)*{\bullet};(37,1)*\txt{6};
(43,4)*{\bullet};(43,1)*\txt{7};
(49,-4);(52,0)**@{-};(55,-4)**@{-};(52,0);(56,4)**@{-};(60,0)**@{-};
(49,-4)*{\bullet};(55,-4)*{\bullet};(60,0)*{\bullet};
(49,-7)*\txt{8};(55,-7)*\txt{9};(60,-3)*\txt{10};
(66,4)*{\bullet};(66,1)*\txt{11};
(72,4)*{\bullet};(72,1)*\txt{12};
\endxy
\\ \\
G:\quad
&\xy
(0,-4);(3,0)**@{-};(6,-4)**@{-};(3,0);(7,4)**@{-};(11,0)**@{-};
(8,-4)**@{-};(11,0);(14,-4)**@{-};
(7,4)*{\bullet};(7,7)*\txt{0};
(20,4)*{\bullet};(20,7)*\txt{1};
(26,4)*{\bullet};(26,7)*\txt{2};
(32,4)*{\bullet};(32,7)*\txt{3};
(38,0);(42,4)**@{-};(46,0)**@{-};(43,-4)**@{-};
(46,0);(49,-4)**@{-};(47,-8)**@{-};(49,-4);(51,-8)**@{-};
(42,7)*\txt{4};(42,4)*{\bullet};
(57,4)*{\bullet};(57,7)*\txt{5};
(63,4)*{\bullet};(63,7)*\txt{6};
(69,4)*{\bullet};(69,7)*\txt{7};
(75,4)*{\bullet};(75,7)*\txt{8};
(81,4)*{\bullet};(81,7)*\txt{9};
(87,4)*{\bullet};(87,7)*\txt{10};
(93,4)*{\bullet};(93,7)*\txt{11};
(99,0);(103,4)**@{-};(107,0)**@{-};
(103,7)*\txt{12};(103,4)*{\bullet};
\endxy
\\ \\
FG:\quad
&\xy
(0,-8);(2,-4)**@{-};(4,-8)**@{-};(2,-4);(5,0)**@{-};(8,-4)**@{-};
(5,0);(9,4)**@{-};(13,0)**@{-};(10,-4)**@{-};(13,0);(16,-4)**@{-};
(6,-8);(8,-4)**@{-};(10,-8)**@{-};
(9,4)*{\bullet};(9,7)*\txt{0};
(22,4)*{\bullet};(22,7)*\txt{1};
(28,-4);(31,0)**@{-};(34,-4)**@{-};
(32,-8);(34,-4)**@{-};(36,-8)**@{-};
(35,-12);(36,-8)**@{-};(37,-12)**@{-};
(31,0);(35,4)**@{-};(39,0)**@{-};
(35,4)*{\bullet};(35,7)*\txt{2};
(45,4)*{\bullet};(45,7)*\txt{3};
(51,4)*{\bullet};(51,7)*\txt{4};
(57,-4);(60,0)**@{-};(63,-4)**@{-};
(60,0);(64,4)**@{-};(68,0)**@{-};
(64,4)*{\bullet};(64,7)*\txt{5};
(74,4)*{\bullet};(74,7)*\txt{6};
(80,0);(84,4)**@{-};(88,0)**@{-};
(84,4)*{\bullet};(84,7)*\txt{7};
\endxy
\end{split}\]

We leave it to the reader to verify that this product is associative
and that the trivial forest is both a left and right identity.  Thus
finite forests form a monoid under this operation.  We extend the
meaning of the symbol \(\scr{F}\) to include this product.

A triple of vertices \((u, u0, u1)\) is called a {\itshape caret}
and is pictured here: \(\xy (-2,-1);(0,1)**@{-};(2,-1)**@{-}\endxy\).
A non-trivial tree is a union of carets.  We add the trivial tree to
the discussion by describing it as the unique tree with zero carets.

Since every finite tree is a union of carets, we can describe a
finite forest as a finite union of carets.  From this it is clear
that the monoid \(\scr{F}\) is generated by the single caret
forests.  Let \(\lambda_i\) be the unique forest with one caret
whose only non-trivial tree is the \(i\)-th tree (which consequently
has only one caret).

From the definition of the product of forests, it is clear that
\(F\lambda_i\) is exactly \(F\) with an extra caret hung from the
\(i\)-the leaf of \(F\).  From this, the following is obvious.

\begin{lemma}\mylabel{ForGen} The forests \(\{\lambda_i\mid i\ge
0\}\) form a generating set for \(\scr{F}\).  \end{lemma}

We let the reader verify the next statement.

\begin{lemma}\mylabel{ForRel} The generators \(\{\lambda_i \mid
i\ge0\}\) satisfy the relations
\(\lambda_q\lambda_m=\lambda_m\lambda_{q+1}\) whenever \(m<q\).
\end{lemma}

To argue that the generating set and relation set of the last two
lemmas form a presentation for \(\scr{F}\), we replace the relation
\(\lambda_q\lambda_m=\lambda_m\lambda_{q+1}\) by the rewriting rule
\mymargin{TheRewriteRule} \begin{equation}\label{TheRewriteRule}
\lambda_q\lambda_m \rightarrow \lambda_m\lambda_{q+1} \quad
\mathrm{whenever}\quad m<q.\end{equation} It is a pleasant exercise
to show that the relation is terminating and locally confluent and
thus complete.  The words that are reduced with respect to
\(\rightarrow\) are the words \(\lambda_{i_0}\lambda_{i_1}\ldots
\lambda_{i_k}\) for which \(i_0\le i_1\le \cdots \le i_k\).  We will
say that words reduced with respect to \tref{TheRewriteRule} are in
normal form.  If \(w\) and \(u\) are two different words in normal
form, then by looking at the leftmost position where they differ, it
is easy to argue that they correspond to two different forests.
This proves the following.

\begin{prop}\mylabel{MonoidFPres} Each element of \(\scr{F}\) is
represented uniquely by a word in the form
\(\lambda_{i_0}\lambda_{i_1}\ldots \lambda_{i_k}\) for which
\(i_0\le i_1\le \cdots \le i_k\).  The monoid \(\scr{F}\) is
presented by \[\langle\lambda_0, \lambda_1, \ldots \mid
\lambda_q\lambda_m = \lambda_m \lambda_{q+1}
\mathrm{\,\,\,whenever\,\,\,} m<q\rangle.\] \end{prop}

In the lemma below, we claim that \(\scr{F}\) is cancellative.
This implies that the equation \(XA=B\) has a unique solution if it
has one at all.  In the case that the equation has a solution, we
write it as \(X=A\backslash B\).  We say that forests \(F\) and
\(G\) are disjoint, if for each \(i\in\N\), at least one of \(F_i\)
or \(G_i\) is trivial.

\begin{lemma}\mylabel{ForPropsLem} The following are true.

{\begin{enumerate}
\renewcommand{\theenumi}{\Roman{enumi}}
\item\mylabel{ForCRM} The monoid \(\scr{F}\) has common
right multiples.  

\item\mylabel{ForCanc} The monoid \(\scr{F}\) is
cancellative.  

\item\mylabel{ForLen} The number of generators that compose
to a given element is a well defined length function on the monoid
\(\scr{F}\).

\item\mylabel{ForTrivUnits} The monoid \(\scr{F}\) has only
trivial units.  

\item\mylabel{ForLCLM} The monoid \(\scr{F}\) has greatest
common right factors, and thus also has least common left multiples.

\item\mylabel{ForGCLF} The monoid \(\scr{F}\) has greatest
common left factors and the greatest common left factor of \(F\) and
\(G\) is \(F\cap G\).  

\item\mylabel{ForLCLMStruct} Let \(F\) and \(G\) be forests
with a common left multiple, and let \(L=PF=QG\) be the least common
left multiple.  The the following are true.  {\begin{enumerate} 
\renewcommand{\theenumii}{\alph{enumii}}
\item \(L\)
is the only least common left multiple of \(F\) and \(G\).  \item If
\(AF=BG\) is a common left multiple of \(F\) and \(G\), then
\(P=(A\cap B)\backslash A\) and \(Q=(A\cap B)\backslash B\).  \item
The forests \(P\) and \(Q\) are disjoint.  \item Each tree in \(L\)
is equal either to a single tree from the forest \(F\) hung on a
trivial tree from \(P\), or to a single tree from the forest \(G\)
hung on a trivial tree from \(Q\).  \end{enumerate}} 

\end{enumerate}}
\end{lemma}

\begin{proof} Mostly left to the reader, and the method of proof can
be done to the reader's taste.  All statements are geometrically
clear from the structure of forests and the nature of the
mutliplication, and they can also be given algebraic proofs from the
relations in lemma \ref{ForRel}.  For example \tref{ForCRM} can be
proven by noting that \(F\cup G\) is a common right multiple for
\(F\) and \(G\), or it can be given an elegant inductive algebraic
proof using the relations of Lemma \ref{ForRel}.  An algebraic proof
for \tref{ForLCLM} can be built by defining a relation on pairs in
\(\scr{F}\) by \((x,y)\rightarrow(z,w)\) if there is a \(\lambda_i\)
so that \(x=z\lambda_i\) and \(y=w\lambda_i\) and showing that it is
complete.  Any common right factor of \(x\) and \(y\) can be built
from a chain from \((x,y)\) to the unique irreducible in the class
containing \((x,y)\).  \end{proof}

\section{The monoid of hedges}

Structures derived from \(\scr{F}\) will use functions defined on
\(\scr{F}\) that factor through a quotient of \(\scr{F}\).  It will
be convenient to be familiar with that quotient.

We said in the introduction that trees will keep track of the order
of splitting of a strand.  If we do not keep track of the order,
then the data in a tree is reduced to a ``shrub.''  A sequence of
shrubs is a hedge and we are about to define the monoid of hedges.

There are many equivalent definitions of a hedge and they each have
their own advantages and disadvantages.  We are less interested in
the details of hedges then we are in their relation to forests, and
we will make all definitions by referring to forests.

We start with the definition that makes the product clear.
Unfortunately, we will rarely refer to this definition in spite of
its advantages.  Let \(F\) be a forest, and let \(l_F:\N\rightarrow
\N\) be defined by \(l_F(i)=j\) if the \(i\)-th leaf of \(F\) is in
\(F_j\).  The function \(l_F\) is a surjection from \(\N\) to
\(\N\), each preimage is finite and non-empty, and only finitely
many preimages have more than one element.  Further \(l_F\)
preserves \(\le\) on \(\N\).  We get a monoid from the set of all
such functions under composition.  It is clear that \(F\mapsto l_F\)
is an epimorphism.  We call \(l_F\) the ``leaf-root'' function of
\(F\).

For the next definition, we note that it is clear that the function
\(l_F\) is completely determined by knowing the size of each set
\(l_F^{-1}(i)\).  This gives a sequence of positive integers, only
finitely many of which are greater than one.  This sequence is just
the sequence \(c_F\) for which \(c_F(j)\) is the number of leaves of
\(F_j\).  Let \(\scr{H}\) be the set of such sequences.  A formula
can be worked out for the product to make
\(\scr{H}\) a monoid isomorphic the the monoid in the previous
paragraph.  In this definition the sequence simply gives the number
of leaves of each ``shrub.''  The epimorphism \(F\mapsto c_F\) is
the ``leaf count'' epimorpihsm from \(\scr{F}\) to \(\scr{H}\).  In
spite of the less pleasant product on \(\scr{H}\), we will refer to
it more often than the others and will use the word ``hedge'' to
refer to an element of \(\scr{H}\).

The third definition takes more information from \(l_F\) and notes
that \(l_F\) is determined by the sets \(l_F^{-1}(i)\).  This is a
partition of \(\N\) into sets each of which is finite and an
interval under \(\le\) on \(\N\).  Further only finitely many sets
have more than one element.  Let \(\scr{P}\) be the set of such
partitions.  Again, a  formula can be worked out
for the product to make \(\scr{P}\) a monoid isomorphic to the
previous two.  This definition has \(l_F^{-1}(i)\) the set of leaf
numbers in \(F_i\).

The proofs of the lemmas in this section are left as exercises for
the reader.

In \(\scr{H}\), we define the hedge \(\nu_i\) by setting
\(\nu_i(i)=2\) and all other values \(1\).

\begin{lemma}\mylabel{HedgeGen} The hedges \(\{\nu_i\mid i\ge
0\}\) form a generating set for \(\scr{H}\).  \end{lemma}

\begin{lemma}\mylabel{HedgeRel} The generators \(\{\nu_i \mid
i\ge0\}\) satisfy the relations \[\nu_q\nu_m = \nu_m\nu_{q+1}
\quad\mathrm{when\,\,\,} m\le q. \] \end{lemma}

From the relations in Lemma \ref{HedgeRel} we derive the rewriting
rules \mymargin{HedgeRewriteRules}
\begin{equation}\label{HedgeRewriteRules}\nu_q\nu_m \rightarrow
\nu_m\nu_{q+1} \quad\mathrm{when\,\,\,} m\le q.\end{equation}

\begin{lemma} The rewriting rules \tref{HedgeRewriteRules} are
locally confluent and terminating and thus complete.  In addition
the inverse rules are also locally confluent and terminating and
thus complete.  \end{lemma}

We end up with two normal forms.  The irreducibles under
\tref{HedgeRewriteRules} are easy to identify as those words
\(\nu_{i_1}\nu_{i_2}\cdots \nu_{i_k}\) with \(i_1<i_2<\cdots <
i_k\).  It is now easy to show that two different irreducible words
represent different hedges.  Thus we have a presentation.

\begin{prop}\mylabel{MonoidHPres} The monoid \(\scr{H}\) is
presented by \[\langle\nu_0, \nu_1, \ldots \mid \nu_q\nu_m = \nu_m
\nu_{q+1} \mathrm{\,\,\,whenever\,\,\,} m\le q\rangle.\] \end{prop}

\begin{lemma}\mylabel{ForHedgeHom} Taking a forest \(F\) to the leaf
count function \(c_F:\N\rightarrow \N-\{0\}\) gives the homomorphism
from \(\scr{F}\) onto \(\scr{H}\) that takes each \(\lambda_i\) to
\(\nu_i\).  \end{lemma}

The irreducibles under the inverse of \tref{HedgeRewriteRules} are
those words \(\nu_{i_1}\nu_{i_2}\cdots \nu_{i_k}\) with \(i_1\ge
i_2\ge \cdots \ge i_k\).  We will call the normal form obtained from
\tref{HedgeRewriteRules} the {\itshape ascending normal form}, and
the normal form obtained from the inverse of
\tref{HedgeRewriteRules} the {\itshape descending normal form}.  It
is more compact to write the descending normal form as
\(\nu_{i_1}^{n_1}\nu_{i_2}^{n_2}\cdots \nu_{i_k}^{n_k}\) where
\(i_1>i_2>\cdots >i_k\) and all \(n_j\) are at least one.

It is a triviality to relate the descending normal form to the
structure of the hedge as a function from \(\N\) to \(\N-\{0\}\).

\begin{lemma}\mylabel{DescendPart} Let the hedge
\(H=\nu_{i_1}^{n_1}\nu_{i_2}^{n_2}\cdots \nu_{i_k}^{n_k}\) be in
descending normal form.  Then \(H(k)=n_j+1\) if \(k=i_j\) for some
\(j\) and \(H(k)=1\) otherwise.  \end{lemma}

\begin{lemma} The monoid \(\scr{H}\) is right cancellative.
\end{lemma}

Since \(\nu_q\ne \nu_{q+1}\) and \(\nu_q\nu_q=\nu_q\nu_{q+1}\), left
cancellativity fails in \(\scr{H}\).

There is a natural isomorphism from the monoid \(\scr{H}\)
consisting of sequences in \(\N-\{0\}\) to the monoid \(\scr{P}\)
consisting of partitions of \(\N\).  The next discussion uses the
advantages of each monoid and we let \(H_P\) be the image in
\(\scr{P}\) of \(H\in \scr{H}\) under the isomorphism.

We define some relations and let the reader verify some claims.  If
\(H\) and \(K\) are hedges, we write \(H\le K\) if for each
\(i\in\N\), we have \(H(i)\le K(i)\).  If \(P\) and \(Q\) are
partitions in \(\scr{P}\), then we write \(P\le Q\) if each set in
\(P\) is contained in some set in \(Q\).

\begin{lemma}\mylabel{FactorFacts} For hedges \(H\) and \(K\) the
following hold.  {\claimenum \item \(H\le K\) if and only if \(H\)
is a left factor if \(K\) (equivalently, \(K\) is a right multiple
of \(H\)).  \item \(H_P\le K_P\) if and only if \(H\) is a right
factor of \(K\) (equivalently, \(K\) is a left multiple of \(H\).
\claimenumend} \end{lemma}

If \(H\) and \(K\) are hedges, then \(\max(H,K)\) is the hedge whose
value at \(i\) is \(\max(H(i), K(i))\) and \(\min(H,K)\) is the
hedge whose value at \(i\) is \(\min(H(i), K(i))\).  If \(P\) and
\(Q\) are partitions in \(\scr{P}\), then there are equivalence
relations \(\sim_P\) and \(\sim_Q\) whose equivalence classes are
\(P\) and \(Q\), respectively.  We set \(P\vee Q\) to be the
partition of classes given by the equivalence relation generated by
\(\sim_P\) and \(\sim_Q\) (that is, by \(\sim_P\cup \sim_Q\)).  We
set \(P\wedge Q\) to be the partition of classes given by the
equivalence relation \(\sim_P\cap \sim_Q\).  Note that \(P\wedge Q\)
and \(P\vee Q\) must be in \(\scr{P}\).

\begin{lemma}\mylabel{ExtremeFacts} Let \(H\) and \(K\) be hedges.
{\claimenum \item The greatest common left factor of \(H\) and \(K\)
is \(\min(H,K)\).  \item The least common right multiple of \(H\)
and \(K\) is \(\max(H,K)\).  \item The greatest common right factor
of \(H\) and \(K\) is the hedge corresponding to \(H_P\wedge K_P\).
\item The least common left multiple of \(H\) and \(K\) is the hedge
corresponding to \(H_P\vee K_P\).  \claimenumend} \end{lemma}

\begin{lemma}\mylabel{PreservLCLM} The homomorphism of Lemma
\ref{ForHedgeHom} takes least common left multiples to least common
left multiples.  \end{lemma}

\section{Incorporating permutations and
braids}\mylabel{BraidingMonoidSec}

We wish to create \ZS products of forests or hedges with braids or
permutations.  The braid group
\(B_n\) on \(n\) strands is as discussed in \cite{birman:braids}.
Since we number things from 0, the strands in \(B_n\) are numbered
0, 1, \dots, \(n-1\), and the generators of \(B_n\) are
\(\sigma_0\), \dots, \(\sigma_{n-2}\).  The infinite braid group
\(B_\infty\) is the direct limit of the \(B_n\) where \(B_n\)
injects into \(B_{n+1}\) by adding a trivial strand at position
\(n\).  The presentation of \(B_\infty\) as a group has generating set
\[\Sigma=\{\sigma_0, \sigma_1, \ldots\}\] 
and relations \mymargin{PermRelsB--C}
\begin{alignat}{2}\label{PermRelsB} \sigma_m\sigma_n &=
\sigma_n\sigma_m, &\qquad&|m-n|\ge2, \\ \label{PermRelsC}
\sigma_m\sigma_{m+1}\sigma_m&=\sigma_{m+1}\sigma_m\sigma_{m+1},
&&m\ge0. \end{alignat}  
The monoid presentation of \(B_\infty\) has
generating set to
\(\Sigma\cup\overline{\Sigma}\) where
\[\overline{\Sigma}=\{\sigma^{-1}_0,\sigma^{-1}_1,\ldots\}\] are the
formal inverses of the elements of \(\Sigma\) and we need add the
relations \mymargin{BraidRelsD} \begin{equation} \label{BraidRelsD}
\sigma_m\sigma^{-1}_m =\sigma^{-1}_m\sigma_m = 1, \qquad
m\ge0. \end{equation}
We follow the convention of \cite{birman:braids} in drawing 
crossings as the following picture of 
\(\sigma_0\sigma_2^{-1}\) shows.
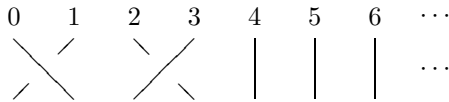
\begin{equation*}
\xy
(0,4);(8,-4)**@{-};
(0,-4);(2,-2)**@{-};
(6,2);(8,4)**@{-};
(16,-4);(24,4)**@{-};
(16,4);(18,2)**@{-};
(22,-2);(24,-4)**@{-};
(32,-4);(32,4)**@{-};
(40,-4);(40,4)**@{-};
(48,-4);(48,4)**@{-};
(56,0)*\txt{\(\cdots\)};
(0,7)*\txt{0};
(8,7)*\txt{1};
(16,7)*\txt{2};
(24,7)*\txt{3};
(32,7)*\txt{4};
(40,7)*\txt{5};
(48,7)*\txt{6};
(56,7)*\txt{\(\cdots\)}
\endxy
\end{equation*}
Also as in \cite{birman:braids}, reading a word in
\(\Sigma\cup\overline\Sigma\) from left to right corresponds to
reading a braid diagram from top to bottom.

The infinite symmetric group \(S_\infty\) is the direct limit of the
finite symmetric groups \(S_n\) and the presentation of \(S_\infty\)
has the generators and relations of \(B_\infty\) in addition to the
relations \mymargin{PermRelsA}\begin{equation}\label{PermRelsA}
\sigma_m^2=1, \qquad m\ge0. \end{equation} Sending each \(\sigma_m\)
in \(B_n\) or \(B_\infty\) to the generator of the same name in
\(S_n\) or \(S_\infty\) gives the standard surjections from braid
groups to symmetric groups.  If \(\sigma\) is in \(B_\infty\), then
the notation \(\sigma(j)\), will always refer to the image of
\(\sigma\) in \(S_\infty\) under this surjection and will give the
image of \(j\) under the permutation.  Context will determine
whether \(\sigma_m\) is a generator of \(B_\infty\) or \(S_\infty\).
The effort it takes to keep track of context will be worth it since
we will be able to deal with both braids and permutations
simultaneously by always using the generating set \(\Sigma\cup
\overline{\Sigma}\).

We regard the \(i\)-th strand of a braid as an arc in 3-space with
top at \((i,0,1)\) and bottom at some \((j,0,0)\).  We say that this
strand has top at \(i\) and bottom at \(j\).  With the usual
interpretation of \(\sigma_m\) in \(S_\infty\) as the transposition
\(m\leftrightarrow m+1\), we get that the \(i\)-th strand of a braid
\(\sigma\) has bottom at \(\sigma(i)\).

\subsection{\protect\ZS products}

To define a \ZS product \(\FB\), we need to put a multiplication on
\(\scr{F}\times B_\infty\) where a generic element will be a forest
followed by a braid.  Following the convention that turns left to
right in word order into top to bottom in a picture, we think of the
braid as hanging from the leaves of the forest.  We build a \ZS
product by telling how \(\beta F\) should be replaced by
\(F'\beta'\) with \(\beta\) and \(\beta'\) from \(B_\infty\) and
\(F\) and \(F'\) from \(\scr{F}\).  The following pictures motivate
the relations we will write down.

\mymargin{ZSRelsA--B}\begin{align}\label{ZSRelsA}
\xy
(0,-7);
(0,-5)**@{-}; (2,-3)**@{-}; (2,5)**@{-};
(4,-7); (4,-5)**@{-}; (2,-3)**@{-};
(6,0)*{\cdots};
(8,-7);
(8,-3)**@{-};(10,-1)**@{-}; (12,1); (14,3)**@{-}; (14,5)**@{-};
(8,5); (8,3)**@{-}; (14,-3)**@{-}; (14,-7)**@{-};
\endxy
\quad\rightarrow\quad
\xy
(0,-7); (0,1)**@{-}; (2,3)**@{-}; (4,1)**@{-}; (4,-7)**@{-};
(2,3);(2,5)**@{-};
(8,-2)*{\cdots};
(10,-7); (10,-5)**@{-};(12,-3)**@{-}; (14,-1); (16,1)**@{-}; (16,5)**@{-};
(10,5); (10,1)**@{-}; (16,-5)**@{-}; (16,-7)**@{-};
\endxy
\qquad\quad&\quad\qquad
\xy
(0,-7);
(0,-5)**@{-}; (-2,-3)**@{-}; (-2,5)**@{-};
(-4,-7); (-4,-5)**@{-}; (-2,-3)**@{-};
(-6,0)*{\cdots};
(-8,-7);
(-8,-3)**@{-};(-10,-1)**@{-}; (-12,1); (-14,3)**@{-}; (-14,5)**@{-};
(-8,5); (-8,3)**@{-}; (-14,-3)**@{-}; (-14,-7)**@{-};
\endxy
\quad\rightarrow\quad
\xy
(0,-7); (0,1)**@{-}; (-2,3)**@{-}; (-4,1)**@{-}; (-4,-7)**@{-};
(-2,3);(-2,5)**@{-};
(-8,-2)*{\cdots};
(-10,-7); (-10,-5)**@{-};(-12,-3)**@{-}; (-14,-1); (-16,1)**@{-};
(-16,5)**@{-}; 
(-10,5); (-10,1)**@{-}; (-16,-5)**@{-}; (-16,-7)**@{-};
\endxy
\\ \notag\\
\label{ZSRelsB}
\xy
(0,-7);(0,-5)**@{-}; (2,-3)**@{-}; (4,-5)**@{-}; (4,-7)**@{-};
(2,-3); (4,-1)**@{-};
(6,1); (8,3)**@{-}; (8,5)**@{-};
(2,5); (2,3)**@{-}; (8,-3)**@{-}; (8,-7)**@{-};
\endxy
\quad\rightarrow\quad
\xy
(-2,-7); (-2,-5)**@{-}; (2,-1)**@{-}; (4,1); (6,3)**@{-}; (8,1)**@{-};
(6,-1)**@{-}; (4,-3); (2,-5)**@{-}; (2,-7)**@{-};
(0,5); (0,3)**@{-}; (8,-5)**@{-}; (8,-7)**@{-};
(6,3);(6,5)**@{-};
\endxy
\qquad\qquad&\qquad\qquad
\xy
(0,-7);(0,-5)**@{-}; (-2,-3)**@{-}; (-4,-5)**@{-}; (-4,-7)**@{-};
(-2,-3); (-4,-1)**@{-};
(-6,1); (-8,3)**@{-}; (-8,5)**@{-};
(-2,5); (-2,3)**@{-}; (-8,-3)**@{-}; (-8,-7)**@{-};
\endxy
\quad\rightarrow\quad
\xy
(2,-7); (2,-5)**@{-}; (-2,-1)**@{-}; (-4,1); (-6,3)**@{-}; (-8,1)**@{-};
(-6,-1)**@{-}; (-4,-3); (-2,-5)**@{-}; (-2,-7)**@{-};
(0,5); (0,3)**@{-}; (-8,-5)**@{-}; (-8,-7)**@{-};
(-6,3);(-6,5)**@{-};
\endxy
\end{align}
Similar pictures motivate relations needed for \(\HB\) and \(\HP\).

To define \(\FB\) and \(\FP\), we let 
\(\Lambda = \{\lambda_0, \lambda_1, \ldots\}\), and to define \(\HB\)
and \(\HP\), we let  \(N=\{\nu_0,
\nu_1,\ldots\}\).  The products \(\FB\) and \(\FP\) 
will be specified by functions
\((\Sigma\cup\overline{\Sigma})\times \Lambda\rightarrow \Lambda\)
written \((\sigma, \lambda)\mapsto \sigma\cdot \lambda\) and
\((\Sigma\cup\overline{\Sigma})\times \Lambda\rightarrow
(\Sigma\cup\overline{\Sigma})\) written \((\sigma, \lambda)\mapsto
\sigma^\lambda\).  Products with \(\scr{H}\) will be specified by
similar functions with \(\Lambda\) replaced by \(N\) and \(\lambda\)
replaced by \(\nu\).  These functions are defined by the following
where \(\epsilon\) represents either \(+1\) or \(-1\):
\mymargin{ZapOpI--III} \begin{align} \label{ZapOpI}
\sigma^\epsilon_q \cdot \lambda_m &= \lambda_{\sigma_q(m)}, \\
\label{ZapOpII} \sigma^\epsilon_q \cdot \nu_m &= \nu_{\sigma_q(m)},
\\ \label{ZapOpIII} (\sigma_q^\epsilon)^{\nu_m} =
(\sigma_q^\epsilon)^{\lambda_m} &= \begin{cases}
\sigma^\epsilon_{q+1}, & m<q, \\ \sigma^\epsilon_q
\sigma^\epsilon_{q+1}, & m=q, \\
\sigma^\epsilon_{q+1}\sigma^\epsilon_q, & m=q+1, \\
\sigma^\epsilon_q, & m>q+1. \end{cases} \end{align}

The reader can check that the relations \(\sigma_q\lambda_m =
\left(\sigma_q\cdot \lambda_m\right)
\left(\sigma_q\right)^{\lambda_m}\) are realizations of the pictures
in \tref{ZSRelsA} and \tref{ZSRelsB}.

We also define a monoid that combines braids with deletion operators
as a \ZS product.  It turns out that deletions from a sequence form a
monoid 
isomorphic to hedges.  We introduce the monoid presentation
\mymargin{MonoidDPres} \begin{equation} \label{MonoidDPres} \scr{Z}=
\langle\delta_0, \delta_1, \ldots \mid \delta_q\delta_m = \delta_m
\delta_{q+1} \mathrm{\,\,\,whenever\,\,\,} m\le q\rangle,
\end{equation} and we let \(\Delta=\{\delta_0, \delta_1,
\ldots\}\).  We will think of \(\delta_q\) as deleting the strand
with top at position \(q\) from a braid, starting at the top.  The
following picture of a spark burning a strand from the top of a
braid motivates the relations we will write down.
\mymargin{SparkPic}\begin{equation}\label{SparkPic}
\xy
(8,8); (8,4)**@{-}; (-4,-8)**@{-}; 
(8,-8); (8,-4)**@{-}; (5,-1)**@{-};
(3,1); (1,3)**@{-}; 
(-1,5); (-4,8)**@{-};
(4,-8); (1,-5)**@{-}; 
(-1,-3); (-4,0)**@{-}; (4,8)**@{-}; (4,8)*{*};
\endxy
\qquad = \qquad
\xy
(8,8); (8,4)**@{-}; (-4,-8)**@{-}; 
(8,-8); (8,-4)**@{-}; (5,-1)**@{-};
(3,1); (-4,8)**@{-};
(4,-8); (1,-5)**@{-}; (-1,-3); (-4,0)**@{-}; (-4,0)*{*};
\endxy
\quad = \quad
\xy
(8,8); (8,4)**@{-}; (-4,-8)**@{-}; 
(8,-8); (8,-4)**@{-}; (5,-1)**@{-};
(3,1); (-4,8)**@{-};
(4,-8)*{*};
\endxy
\end{equation}

Products \(\BZ\) and \(\PZ\) 
will be specified by functions
\(\Delta\times(\Sigma\cup \overline{\Sigma})\rightarrow \Delta\)
written \((\delta,\sigma)\mapsto \delta\cdot \sigma\) and
\(\Delta\times(\Sigma\cup \overline{\Sigma})\rightarrow (\Sigma\cup
\overline{\Sigma})\) written \((\delta,\sigma)\mapsto
\delta^\sigma\).  These functions are defined by the following where
\(\epsilon\) represents either \(+1\) or \(-1\):
\mymargin{ZapOpIV--V} \begin{align} \label{ZapOpIV}
(\delta_q)^{\sigma_m^\epsilon} &= \delta_{\sigma_m(q)}, \\
\label{ZapOpV} \delta_q\cdot \sigma_m^\epsilon &=
\begin{cases}\sigma^\epsilon_{m-1}, \quad &q<m, \\ 1, &q=m,m+1, \\
\sigma^\epsilon_m, & q>m+1. \end{cases} \end{align}
The illustration in \tref{SparkPic} shows the truth of
\(\delta_1\sigma_0\sigma_1\sigma_0 = \delta_0\sigma_1\sigma_0 =
\sigma_0\delta_0\sigma_0 = \sigma_0\delta_1\).  The \(\delta_1\) at
the end is to be interpreted as ready to delete strand 1 from any
braid that might be multiplied on the right of the original.

In \tref{ZapOpIII} and \tref{ZapOpV}, the consistent treatment of
the exponent \(\epsilon\) allows restriction of domain and codomain
from \(\Sigma\cup\overline{\Sigma}\) to \(\Sigma\) when working with
\(S_\infty\) instead of \(B_\infty\).  

The formula \tref{ZapOpIV} has to be interpreted carefully.  It is
only a statement about expressions involving generators.  It is not
meant to imply that it applies to words in these generators.  In
fact, if \tref{ZapOpIV} is followed literally, then we get
\mymargin{ActnChkI}\begin{equation}\label{ActnChkI}
((\delta_q)^{\sigma_m})^{\sigma_n} =
(\delta_{\sigma_m(q)})^{\sigma_n} = \delta_{\sigma_n(\sigma_m(q))}
\end{equation} which is to be compared with the incorrect
\mymargin{ActnChkII} \begin{equation}\label{ActnChkII}
(\delta_q)^{(\sigma_m\sigma_n)} = \delta_{(\sigma_m\sigma_n)(q)}.
\end{equation} 

From Lemma \ref{ZapPresIntCompl} we know we get a consistent action
of words in the \(\sigma\) on the various \(\delta\) if we define
the left side of \tref{ActnChkII} to equal the right side of
\tref{ActnChkI}.  The various elements \(\sigma\) as show up in the
subscripts in \tref{ActnChkI} are to be interpreted as permutations
(and are transpositions) and the reverse of a string of
transpositions is the inverse of the original string.  This gives
one conclusions of the following lemma.  The other conclusions
follow from Lemma \ref{ZapPresIntCompl} and from \tref{ZapOpI} and
\tref{ZapOpII} without the complications relating to \tref{ZapOpIV}
since in \tref{ZapOpI} and \tref{ZapOpII}, the various \(\sigma\)
act on the left.

\begin{lemma}\mylabel{HowSigsAct} For any \(\lambda_j\in \Lambda\),
\(\nu_j\in N\) or \(\delta_j\in \Delta\) and \(\tau\) in
\(S_\infty\) or \(B_\infty\), we have \[\tau\cdot \lambda_j =
\lambda_{\tau(j)}, \qquad \tau\cdot \nu_j = \nu_{\tau(j)}, \qquad
(\delta_j)^\tau = \delta_{\tau^{-1}(j)}.\] \end{lemma}

\begin{prop}\mylabel{AreZappaProds} The functions defined by
\tref{ZapOpI}, \tref{ZapOpII}, \tref{ZapOpIII}, \tref{ZapOpIV} and
\tref{ZapOpV} define Zappa products \(\scr{F}\zapprod S_\infty\),
\(\scr{F}\zapprod B_\infty\), \(\scr{H}\zapprod S_\infty\),
\(\scr{H}\zapprod B_\infty\), \(S_\infty\zapprod\scr{Z}\) and
\(B_\infty \zapprod \scr{Z}\).  \end{prop}

\begin{proof} This is an orgy of checking the requirements of Lemma
\ref{TwistedPresIII} which is left to the reader.
The number of cases is large.  Note that the roles of
\(\langle X\mid R\rangle\) and \(\langle Y\mid T\rangle\) in that
lemma must be reversed in dealing with the products with
\(\scr{Z}\).  We point out that 
the flexibility of the hypothesis of Lemma \ref{TwistedPresIII} that
allows either of \((\sigma\cdot u, \sigma\cdot v)\) or
\((\sigma\cdot v, \sigma\cdot u)\) to be a relation for \(\scr{F}\)
if \((u,v)\) is a relation for \(\scr{F}\) must be used when showing
that the related pair
\((\lambda_{m+1}\lambda_m,\lambda_m\lambda_{m+2})\) in \(\scr{F}\)
is carried to the related pair \((\sigma_m\cdot
(\lambda_{m+1}\lambda_m), \sigma_m\cdot (\lambda_m\lambda_{m+2})) =
(\lambda_m\lambda_{m+2}, \lambda_{m+1}\lambda_m).\) \end{proof}

There are pictures that can go with the products of Proposition
\ref{AreZappaProds}.  For example, elements of \(\scr{F}\zapprod
B_\infty\) can be thought as forests with braids hanging from the
leaves.  We show a calculation (only the first tree of each forest is
shown) 
\[
\left(\,\,\xy 
(4,-6); (-4,2)**@{-}; (0,6)**@{-}; (4,2)**@{-}; 
(1,-1)**@{-}; (-1,-3); (-4,-6)**@{-};
\endxy\,\,\right)^2
\qquad = \qquad
\xy
(4,-12); (4,0)**@{-}; (-4,8)**@{-}; (0,12)**@{-}; (4,8)**@{-}; (1,5)**@{-};
(-1,3); (-4,0)**@{-}; 
(0,-12); (-8,-4)**@{-}; (-4,0)**@{-}; (0,-4)**@{-}; (-3,-7)**@{-};
(-5,-9); (-8,-12)**@{-};
\endxy
\qquad = \qquad
\xy
(8,-10); (8,-6)**@{-}; (-4,6)**@{-}; (0,10)**@{-}; (8,2)**@{-}; (5,-1)**@{-};
(3,-3); (1,-5)**@{-}; (-1,-7); (-4,-10)**@{-};
(4,6); (1,3)**@{-}; (-1,1); (-4,-2)**@{-}; (4,-10)**@{-};
\endxy
\]
or equivalently, \(\lambda_0\sigma_0\lambda_0\sigma_0 =
\lambda_0\lambda_1 \sigma_0\sigma_1\sigma_0\).

The check that the relation \(\nu_m\nu_m=\nu_m\nu_{m+1}\) cooperates
with \tref{ZapOpIII} shows that the action of forests on braids
successfully factors through the action of hedges on braids.

\subsection{Some algebraic properties of the products} We would like
to prove that \(\FP\) and \(\FB\) share some of the properties that
are possessed by \(\scr{F}\).  We will make use of Lemma
\ref{PartMulProps}, so we start by verifying some of the properties
needed by that lemma.  We first need some technical lemmas.

\begin{lemma}\mylabel{CoActions} The following equalities hold
concerning the actions \(B_\infty\times \scr{F}\rightarrow
B_\infty\) used in creating \(\scr{F}\zapprod B_\infty\) and
\(\scr{Z}\times B_\infty\rightarrow B_\infty\) use in creating
\(B_\infty\zapprod \scr{Z}\): \[\delta_{\sigma_q(m)}\cdot
\left((\sigma_q)^{\lambda_m}\right) = \sigma_q,
\qquad\mathrm{and}\qquad \left(\delta_{\sigma_q(m)}\right) ^{
\left((\sigma_q)^{\lambda_m}\right)}=\delta_m.\] \end{lemma}

\begin{proof} We write out the calculation for \(m=q\).
\[\begin{split} \delta_{\sigma_q(q)}\cdot
\left((\sigma_q)^{\lambda_q}\right) &= \delta_{q+1}\cdot
(\sigma_q\sigma_{q+1}) = (\delta_{q+1}\cdot
\sigma_q)((\delta_{q+1})^{\sigma_q}\cdot \sigma_{q+1}) \\ &=
(1)(\delta_q\cdot \sigma_{q+1}) = \sigma_q. \end{split}\]
\[\begin{split} \left(\delta_{\sigma_q(q)}\right) ^{
\left((\sigma_q)^{\lambda_q}\right)} &=
(\delta_{q+1})^{(\sigma_q\sigma_{q+1})} =
(\delta_{\sigma_q(q+1)})^{\sigma_{q+1}} \\ &=
(\delta_q)^{\sigma_{q+1}} = \delta_{\sigma_{q+1}(q)} = \delta_q =
\delta_q. \end{split}\] The other cases, \(m<q\), \(m=q+1\) and
\(m>q+1\) are left to the reader. \end{proof}

\begin{lemma}\mylabel{CoActExtn} In the setting of Lemma
\ref{CoActions} we have \[\delta_{\tau(m)}\cdot
\left(\tau^{\lambda_m}\right) = \tau\] for any \(\tau\in B_\infty\).
\end{lemma}

\begin{proof} From Lemma \ref{CoActions}, the result holds if
\(\tau\) is a single generator.  Consider \(\tau=\sigma\omega\) for
some generator \(\sigma\) so that \(\omega\) has shorter length then
\(\tau\).  Then \[\begin{split} \delta_{(\sigma\omega)(m)} \left (
(\sigma\omega)^{\lambda_m} \right) &= \delta_{(\sigma\omega)(m)}
\left( \sigma^{\omega\cdot \lambda_m}\omega^{\lambda_m}\right) \\ &=
\left( \delta_{\sigma(\omega(m))} \cdot
\left(\sigma^{\lambda_{\omega(m)}} \right) \right) \left( \left (
\delta_{\sigma(\omega(m))}
\right)^{\left(\sigma^{\lambda_{\omega(m)}} \right)} \cdot
(\omega^{\lambda_m}) \right) \\ &=(\sigma)\left(\delta_{\omega(m)}
\cdot (\omega^{\lambda_m}) \right) \\ &= \sigma\omega=\tau
\end{split}\] where the next to last line is justified by the two
parts of Lemma \ref{CoActions} and the last line is by induction
since the length of \(\omega\) is less than that of \(\tau\).
\end{proof}

\begin{cor}\mylabel{CoActInj} In the setting of Lemma
\ref{CoActions}, for each \(u\in \scr{F}\), the family of functions
\(B_\infty\times \scr{F}\rightarrow B_\infty\) given by
\((\tau,u)\mapsto \tau^u\) is a family of injections.  \end{cor}

\begin{proof}  By Lemma \ref{CoActExtn}, this is true if \(u\) is
some \(\lambda_m\).  The claim follows since
\(\tau^{uv}=(\tau^u)^v\) and a composition of injections is an
injection.  \end{proof}

We are now in a position to prove the following facts about
\(\FB\).

\begin{prop}\mylabel{FBOreProps} The monoid \(\FB\) is cancellative
and has common right multiples.  Further, declaring the length of
\(G\tau\) with \(G\in \scr{F}\) and \(\tau\in B_\infty\) to be the
length of \(G\) as given in Lemma \ref{ForPropsLem}\tref{ForLen}
gives a length function on \(\FB\).  \end{prop}

\begin{proof} For cancellativity, Lemma \ref{PartMulProps} and the
unstated version for left cancellativity says that we need that
\((\tau,u)\mapsto \tau^u\) and \((\tau,u)\mapsto \tau\cdot u\) are
injective families.  We get one from Corollary \ref{CoActInj} and
the other from the fact that \(B_\infty\) is a group and that
\((\tau,u)\mapsto \tau\cdot u\) is an action.  For common right
multiples, Lemma \ref{PartMulProps} requires that \((\tau,u)\mapsto
\tau\cdot u\) is a surjective family.  This follows from the fact
that we have a group action.  That the claimed length function for
\(\FB\) is truly a length function follows from the fact that
\(B_\infty\) is a group whose action on \(\scr{F}\) takes generators
to generators and is thus length preserving.  \end{proof}

Identical arguments over the last few lemmas repeated for
\(S_\infty\) give the following.

\begin{prop}\mylabel{FPOreProps} The monoid \(\FP\) is cancellative
and has common right multiples.  Further, declaring the length of
\(G\tau\) with \(G\in \scr{F}\) and \(\tau\in S_\infty\) to be the
length of \(G\) as given in Lemma \ref{ForPropsLem}\tref{ForLen}
gives a length function on \(\FP\).  \end{prop}

\subsection{Least common left multiples}  Least common left
multiples are needed to get reduced terms in groups of fractions.
Here they require extra work.

We need more information than given in Lemmas \ref{CoActions} and
\ref{CoActExtn}.  The content of Lemma \ref{CoActExtn} is that going
from \(\tau\in B_\infty\) to \(\tau^{\lambda_m}\) splits a strand
into two parallel strands, and \(\delta_{\tau(m)}\) restores
\(\tau\) by deleting one of the parallel strands.  We elaborate on
that.  What follows are discussions about inductive extentions of
Lemmas \ref{CoActions} and \ref{CoActExtn} from statements about the
behavior of generators to statements about the behavior of arbitrary
elements.

If \(\sigma\) is a braid representative, then we say strands \(i\)
and \(i+1\) are {\itshape companions} if \(\sigma(i+1)=\sigma(i)\).
In this case there is a circle \(J\) built from strands \(i\) and
\(i+1\), the straight line segement joining the tops of the strands,
and the straight line segement joining the bottoms of the strands.
If there is a disk that meets each plane \(z=t\), \(0\le t\le 1\),
in a single line segment of length 1 that is parallel to the
\(x\)-axis, whose boundary is \(J\) and which does not meet any
strand of \(\sigma\) in its interior, then we say that strands \(i\)
and \(i+1\) are {\itshape parallel}.

If a braid has a representative in which strands \(i\) and \(i+1\)
are parallel, then we say that strands \(i\) and \(i+1\) are
{\itshape weakly parallel} in any other representative.  It is
possible to characterize weakly parallel strands by defining a
winding number of two strands that are companions and showing that
companion strands are weakly parallel if they have winding number 0
and the circle \(J\) of the previous paragraph bounds a disk whose
interior is disjoint from the strands in \(\sigma\).

If \(\sigma\) is a braid representative, then a partition of the
strands of \(\sigma\) is into {\itshape weak parallel classes} if
any two consecutive braids in a class are weakly parallel.  We
insist that elements of the partition be finite.  Since we consider
braids in \(B_\infty\), we do not insist that these classes be
maximal.  We further insist that only finitely many classes have
more than one strand.

If we label strands by their strand numbers, then a partition of the
strands of a braid is identified with a partition of \(\N\).  The
following is straightforward.

\begin{lemma} If \(\sigma\) is a braid representative, and \(C\) is
a partition of \(\N\) into weak parralel classes of \(\sigma\), then
there is a representative \(\sigma'\) of the same braid in which any
two consecutive strands in the same class of \(C\) are parallel.
\end{lemma}  

We now drop the phrase ``weakly parallel'' and ``weak parallel
classes'' and only refer to parallel strands and parallel classes,
and we think of the property ``parallel'' as being attached to
strands of of a braid not the strands of a representative.

Recall that partitions of \(\N\) into finite sets of consecutive
numbers with all but finitely many sets of size one can be viewed as
hedges.  The next lemma refers to the operations on hedges as used
in Lemma \ref{ExtremeFacts}.

\begin{lemma}\mylabel{CombineParClasses} If partitions \(C\) and
\(D\) of \(\N\) are partitions of the strands of \(\sigma\in
B_\infty\) into parallel classes, then so is the partition \(C\vee
D\).  \end{lemma}

The following is geometrically ``obvious'' and is proven
inductively, first on the number of generators in the braid, and
then on the number of generators in the hedge.  The start of the
induction is from the pictures in \tref{ZSRelsB}.

\begin{lemma}\mylabel{BraidParHedgeCor} Let \(\sigma\) be a braid
and let \(u\) be a hedge corresponding to partition \(Q\).  Then
\(Q\) is a partition of \(\sigma^u\) into parallel classes.
\end{lemma}

If \(\sigma\) is a braid, and \(Q\) a partition into parallel
classes of strands, then we can ``collapse'' each class into a
single strand.  This is accomplished by deleting all strands but one
in each class.  It is clear that it does not matter which strand is
the one in each class chosen to remain.  We use \(\sigma/Q\) to
denote the result of this operation.  The next lemma is again proven
by induction, first on the number of generators of the braid and
then on the number of generators of the hedge corresponding to
\(Q\).

\begin{lemma}\mylabel{CollapseIsInj} Let \(\sigma\) and \(\tau\) be
braids so that a partition \(Q\) into finite sets, only finitely
many of which are not singletons, is a partition of both \(\sigma\)
and \(\tau\) into parallel classes.  If \(\sigma/Q=\tau/Q\),
then \(\sigma=\tau\).  \end{lemma}

The next lemma is built inductively from Lemma \ref{CoActExtn}.

\begin{lemma}\mylabel{CancelHedge} Let \(\sigma\) be a braid and let
\(u\) be a hedge with corresponding partition \(Q\) of \(\N\).  Then
\((\sigma^u)/Q=\sigma\).  \end{lemma}

\begin{lemma}\mylabel{HedgeReconstruct} Let \(\sigma\) be a braid
with a partition \(Q\) into parallel classes with each class
finite and only finitely many classes not singletons.  Let \(u\) be
the hedge corresponding to \(Q\).  Then \((\sigma/Q)^u=\sigma\).
\end{lemma}

\begin{proof} We have \(((\sigma/Q)^u)/Q=\sigma/Q\) by Lemma
\ref{CancelHedge}.  But \((\sigma/Q)^u\) has \(Q\) as a partition
into parallel classes by Lemma \ref{BraidParHedgeCor}.  Now we get
the conclusion from Lemma \ref{CollapseIsInj}.  \end{proof}

\begin{lemma} \mylabel{StrCoCon} The family of functions
\(B_\infty\times \scr{F}\rightarrow B_\infty\) written
\((\sigma,u)\mapsto \sigma^u\) used in creating \(\FB\) is a
strongly coconfluent family of injections.  \end{lemma}

\begin{proof} The injective properties follow from Lemma
\ref{CancelHedge}.  For the coconfluence, we must show that if
\(\sigma^u=\tau^v\) where \(u\) and \(v\) have a common left
multiple, then there is a braid \(\beta\) and \(p\) and \(q\) so
that \(pu=qv\), \(\beta^p=\sigma\) and \(\beta^q=\tau\) all hold.
From Lemma \ref{StrongCoConf}, we know that if this holds, then it
will hold when \(w=pu=qv\) is the least common left multiple of
\(u\) and \(v\), so we assume that it is.  (We know that the least
common left multiple of \(u\) and \(v\) must exist.)  We denote the
homomorphism from \(\scr{F}\) to \(\scr{H}\) by \(u\mapsto
\overline{u}\), and we note that \(\overline{w} =
\overline{p}\overline{u} = \overline{q}\overline{v}\) is the least
common left multiple of \(\overline{u}\) and \(\overline{v}\) in
\(\scr{H}\) by Lemma \ref{PreservLCLM}.

Let \(\alpha=\sigma^u=\tau^v\).  Since the action of \(\scr{F}\)
factors through \(\scr{H}\), we note that
\(\alpha=\sigma^{\overline{u}}=\tau^{\overline{v}}\).  We let
\(C_u\) and \(C_v\) denote the partition of \(\N\) corresponding to
\(\overline{u}\) and \(\overline{v}\), respectively, and we note
that both \(C_u\) and \(C_v\) are partitions of the strands of
\(\alpha\) into weak parallel classes.  Thus \(C_u\vee C_v\) must be
a partition of the strands of \(\alpha\) into weak parallel classes.

From Lemma \ref{ExtremeFacts}, we know that the hedge corresponding
to \(C_u\vee C_v\) is the least common left multiple of
\(\overline{u}\) and \(\overline{v}\) which is \(\overline{w}\).
Let \(C=C_u\vee C_v\).  From Lemma \ref{HedgeReconstruct}, we have
\((\alpha/C)^{\overline{w}} =\alpha\).  Let \(\beta=\alpha/C\).

Now \(\sigma^{\overline{u}} = \alpha = \beta^{\overline{w}} =
\beta^{\overline{p}\overline{u}}\) and by Corollary \ref{CoActInj},
we get \(\sigma=\beta^{\overline{p}} = \beta^p\).  Similarly, we get
\(\tau = \beta^q\).  This completes the proof.  \end{proof}

\begin{prop}\mylabel{FBFPLCLM} The monoids \(\FB\) and \(\FP\) have
least common left multiples.  Further, the least common left
multiple \((p,\alpha)(u,\sigma) = (q,\beta)(v,\tau)\) of
\((u,\sigma)\) and \((v,\tau)\) can be constructed so that
\(p(\alpha\cdot u)=q(\beta\cdot v)\) is the least common left
multiple of \((\alpha\cdot u)\) and \((\beta\cdot v)\) in
\(\scr{F}\).  \end{prop}

\begin{proof} This follows from Lemma \ref{PartMulPropsLCLM} and
from Lemma \ref{StrCoCon} and a corresponding lemma for
\(\FP\) which has an identical proof given the almost identical
behavior of \(S_\infty\) and \(B_\infty\).  \end{proof}

Lemma \ref{TwistedPresIII} gives presentations for 
\(\FP\) and \(\FB\)
as follows:
\mymargin{FPPres} \begin{equation}\label{FPPres} \begin{split}
\FP = \langle \Lambda \cup \Sigma \mid \,& \lambda_q\lambda_m =
\lambda_m \lambda_{q+1},\qquad m<q, \\ &\sigma_m^2=1, \qquad m\ge0,
\\ &\sigma_m\sigma_n = \sigma_n\sigma_m, \qquad |m-n|\ge2, \\
&\sigma_m\sigma_{m+1}\sigma_m = \sigma_{m+1}\sigma_m\sigma_{m+1},
\qquad m\ge0, \\ &\sigma_q\lambda_m = (\sigma_q\cdot
\lambda_m)(\sigma_q)^{\lambda_m}\qquad\rangle,
\end{split}\end{equation} and \mymargin{FBPres}
\begin{equation}\label{FBPres} \begin{split} \FB = \langle
\Lambda \cup \Sigma \cup \overline{\Sigma} \mid \,
&\lambda_q\lambda_m = \lambda_m \lambda_{q+1},\qquad m<q, \\
&\sigma_m\sigma^{-1}_m=1, \qquad m\ge0, \\ &\sigma_m\sigma_n =
\sigma_n\sigma_m, \qquad |m-n|\ge2, \\ &\sigma_m\sigma_{m+1}\sigma_m
= \sigma_{m+1}\sigma_m\sigma_{m+1}, \qquad m\ge0, \\
&\sigma^\epsilon_q\lambda_m = (\sigma^\epsilon_q \cdot \lambda_m)
(\sigma^\epsilon_q)^{\lambda_m}\qquad
\rangle. \end{split}\end{equation}

\section{Groups of fractions}\mylabel{MonoidFracSec}

The monoids \(\FB\) and \(\FP\) are cancellative with common right
multiples and thus have groups of right fractions.  We let
\(\widehat{BV}\) be the group of right fractions for \(\FB\) and we
let \(\widehat{V}\) be the group of right fractions for \(\FP\).

\subsection{Embeddings} The group of fractions construction and the
\ZS product both involve embeddings (Proposition \ref{Ore} and Lemma
\ref{PartMulPropsII}).  This is reflected in the following where we
use notation based on the fact that if \(\scr{M}\) is cancellative
monoid with common right multiples, then elements of the group of
right fractions of \(\scr{M}\) are represented by elements of
\(\scr{M}\times\scr{M}\).

\begin{prop}\mylabel{EmbedsProp} Sending \(F\) to \((F,1)\) embeds
\(\scr{F}\) into \(\FP\) and \(\FB\).  Sending \(\sigma\) to
\((1,\sigma)\) embeds \(S_\infty\) into \(\FP\) and embeds
\(B_\infty\) into \(\FB\).  Sending \((F,\beta)\) to
\(((F,\beta),1)\) embeds \(\FP\) in \(\widehat{V}\) and embeds
\(\FB\) into \(\widehat{BV}\).  \end{prop}

\subsection{Some presentations} It is easy to give infinite
presentations of \(\widehat{V}\) and \(\widehat{BV}\).

From Propositions \ref{Ore} and \ref{MonoidFPres} and from
\tref{ZapOpI}, \tref{ZapOpIII}, \tref{FPPres} and \tref{FBPres}, we
get the following where \(\Lambda = \{\lambda_0, \lambda_1,
\ldots\}\) and \(\Sigma=\{\sigma_0, \sigma_1, \ldots\}\).

\begin{thm}\mylabel{InfinitePres} The groups
\(\widehat{V}\) and \(\widehat{BV}\) are presented as groups by
\begin{alignat*}{2}
\widehat{V} = 
  \langle
    \Lambda \cup \Sigma \mid \,
    & \lambda_q\lambda_m = \lambda_m \lambda_{q+1}, 
      &&m<q, \\ 
    &\sigma_m^2=1, 
      &&m\ge0, \\
    &\sigma_m\sigma_n = \sigma_n\sigma_m, 
      &&|m-n|\ge2, \\
    &\sigma_m\sigma_{m+1}\sigma_m = \sigma_{m+1}\sigma_m\sigma_{m+1},
      &&m\ge0, \\ 
    &\sigma_q\lambda_m = \lambda_m\sigma_{q+1},
      &&m<q, \\
    &\sigma_m\lambda_m = \lambda_{m+1}\sigma_m\sigma_{m+1},
      &&m\ge0, \\
    &\sigma_m\lambda_{m+1} = \lambda_m\sigma_{m+1}\sigma_m,
      &&m\ge0, \\
    &\sigma_q\lambda_m = \lambda_m\sigma_q, 
      &&m>q+1
  \rangle, \\ \\
\widehat{BV} =
  \langle 
    \Lambda \cup \Sigma \mid \, 
    &\lambda_q\lambda_m = \lambda_m \lambda_{q+1}, 
      &&m<q, \\ 
    &\sigma_m\sigma_n = \sigma_n\sigma_m,
      &&|m-n|\ge2, \\ 
    &\sigma_m\sigma_{m+1}\sigma_m = \sigma_{m+1}\sigma_m\sigma_{m+1},
      &&m\ge0, \\
    &\sigma^\epsilon_q\lambda_m = \lambda_m\sigma^\epsilon_{q+1},
      &&m<q,\,\,\epsilon=\pm1, \\
    &\sigma^\epsilon_m\lambda_m = 
          \lambda_{m+1}\sigma^\epsilon_m\sigma^\epsilon_{m+1},
      &&m\ge0,\,\,\epsilon=\pm1, \\
    &\sigma^\epsilon_m\lambda_{m+1} = 
          \lambda_m\sigma^\epsilon_{m+1}\sigma^\epsilon_m,
      &&m\ge0,\,\,\epsilon=\pm1, \\
    &\sigma^\epsilon_q\lambda_m = \lambda_m\sigma^\epsilon_q, 
      &&m>q+1,\,\,\epsilon=\pm1
  \rangle. 
\end{alignat*} 
\end{thm}

Some of the relations are redundant.  The relations
\(\sigma_m\lambda_{m+1} = \lambda_m\sigma_{m+1}\sigma_m\) follow
from the relations \(\sigma_m\lambda_m =
\lambda_{m+1}\sigma_m\sigma_{m+1}\) in \(\widehat{V}\) by bringing
each \(\sigma_m\) and \(\sigma_{m+1}\) to the other side fo the
equality.  Similarly, the relations \(\sigma^\epsilon_m\lambda_{m+1}
= \lambda_m\sigma^\epsilon_{m+1}\sigma^\epsilon_m\) follow from the
relations \(\sigma^\epsilon_m\lambda_m =
\lambda_{m+1}\sigma^\epsilon_m\sigma^\epsilon_{m+1}\) in
\(\widehat{BV}\).  Also, the exponents \(\epsilon\) can be
eliminated from several of the relations in \(\widehat{BV}\) because
of the group setting.

\subsection{Normal forms}

The monoids have least common left multiples and length functions.
Thus elements in \(\widehat{BV}\) and \(\widehat{V}\) have
representatives of the fractions in reduced terms.  The next lemma
gives the details of the normal form that comes from the reduced
terms.  We refer to the length function on the monoid \(\scr{F}\)
from Lemma \ref{ForPropsLem}\tref{ForLen}.

\begin{thm}\mylabel{BigTNormalForms} (I) Each element \(x\) of
\(\widehat{V}\) or \(\widehat{BV}\) is represented
uniquely by a triple \((G, \alpha, H)\) with the conditions that
{\claimenum \item \(G\) and \(H\) are in \(\scr{F}\), \item
\(\alpha\) is in 
\(S_\infty\) if \(x\in \widehat{V}\), and is in \(B_\infty\) if
\(x\in \widehat{BV}\), \item \(x=(G\alpha)(H)^{-1}\), \item the
length of \(G\) is minimal among all triples satisfying (a--c).
\claimenumend} 

(II) Any other representative of \(x\) is of the form
\((G\alpha J\gamma)(H J\gamma)^{-1}\) for some \(J\) in \(\scr{F}\)
and \(\gamma\) in the appropriate one of \(S_\infty\) or
\(B_\infty\).

(III) The triple \((G, \alpha, H)\) of (I) is characterized by the
fact that \(x=G\alpha H^{-1}\) and for no \(G'\), \(H'\) and
\(\lambda_i\) in \(\scr{F}\) and \(\alpha'\) in
\(S_\infty\) or \(B_\infty\), as appropriate, is it true that
\(G\alpha=G'\alpha'\lambda_i\) and \(H=H'\lambda_i\).  \end{thm}

\begin{proof} We can do both groups at once if we consider
expressions such as \(G\mu\), \(H\tau\) or \(J\gamma\) in which
\(G\), \(H\) and \(J\) are in \(\scr{F}\) and \(\mu\), \(\tau\) and
\(\gamma\) are either in \(S_\infty\) or
\(B_\infty\) depending on whether we are discussing, respectively,
\(\widehat{V}\) or \(\widehat{BV}\).  Take an
element \(x\) in one of the groups and represent it as
\(x=(G\mu)(H\tau)^{-1}\) so that the length of \(G\mu\) is minimal
among the representatives of \(x\).  From Lemma
\ref{ForPropsLem}\tref{ForLen} and Propositions
\ref{FBOreProps} and \ref{FPOreProps}, the length of \(G\mu\) is the
length of \(G\) as a word in the symbols \(\lambda_i\).  Thus
\(x=(G\mu\tau^{-1})(H)^{-1}\) is another representative of \(x\)
with the same properties, and we take the desired triple to be \((G,
\alpha, H)\) with \(\alpha=\mu\tau^{-1}\).  This satisfies (a--d)
and We now need to consider uniqueness.

Since the length of \(G\alpha\) is minimal, we know from Lemma
\ref{ReducedTerms}, that \((G\alpha)(H)^{-1}\) is in reduced terms.
If \((G'\alpha')(H')^{-1}\) is another representative of \(x\)
coming from a triple \((G', \alpha', H')\) satisfying (a--d), then
there is some \(J\gamma\) so that \(G'\alpha'=G\alpha
J\gamma=(G(\alpha \cdot J))(\alpha^J\gamma)\) and \(H'=H J\gamma\).
From \tref{ZapOpI}, the action of \(\alpha\) on \(J\) preserves the
length of \(J\), so the minimality of the lengths of \(G\) and
\(G'\) force the length of \(J\) to be zero.  The only element of
\(\scr{F}\) with length zero is the identity.  The uniqueness of
representation in a \ZS product forces \(\gamma\) to be the
identity.  This finishes (I).

Item (II) follows from the fact that \((G\alpha)(H)^{-1}\) is in
reduced terms.

To see (III), we note that from (II), the length of \(G\) is minimal
if the test in (III) is satisfied, and the length is not minimal if
the test in (III) is not satsified.  \end{proof}

\subsection{An isomorphism} We argue that each element \(F\beta\) of
\(\FB\) gives an element in the geometric description of
\(\widehat{BV}\) from Section \ref{GeomDescrSec}.  The forest \(F\)
tells how to break the intervals in \(\scr{J}=\{[2i,\,2i+1]\mid
i\ge0\}\) into smaller intervals.  The tree \(F_i\) gives
instructions on breaking the interval \([2i,\,2i+1]\).  The braid
\(\beta\) tells how to reorder the intervals by an isotopy of
\(\R^2\).  The image intervals are now resized and moved
horizontally on the \(x\)-axis until each maps affinely onto one of
the \([2i,\,2i+1]\) so that they are all covered.  Thus \(F\beta\)
can be thought of as a braiding that takes the subdivided intervals
from \(\scr{J}\) onto the unsubdivided intervals from \(\scr{J}\).
The fact that the forest \(F\) is finite corresponds to the fact
that in the description of Section \ref{GeomDescrSec}, all but
finitely many intervals of the cover of \(X\) must be intervals from
\(\scr{J}\).  It is clear that any element of the group from Section
\ref{GeomDescrSec} is of the form \((F\beta)(G\gamma)^{-1}\) for
some pair of elements \(F\beta\) and \(G\gamma\) from \(\FB\).

The multiplication of forests corresponds to successive subdivisions
of intervals and the relations from Lemma \ref{ForRel} on
\(\scr{F}\) are seen to hold among the subdivision operations.  The
braid relations \tref{PermRelsB} and \tref{PermRelsC}  hold for the
braiding and the \ZS relations \tref{ZapOpI} and \tref{ZapOpIII}
hold as pictured in \tref{ZSRelsA} and \tref{ZSRelsB}.  Thus the
association of elements of \(\FB\) to braidings of intervals is a
homomorphism.  By Proposition \ref{Ore}(d), the homomorphism extends
to one defined on \(\widehat{BV}\).  As remarked in the previous
paragraph, the homomorphism is a surjection.

If \((F,\beta, G)\) from \(\widehat{BV}\) is taken to the identity
in the group of Section \ref{GeomDescrSec}, then we have a bijection
of interval collections that must be the identity.  It is easy to
argue that different forests give different collections of
intervals, so \(F=G\).  By conjugating by \(F\), we see that
\((1,\beta,1)\) is also taken to the identity.  But this is just a
braiding of the intervals in \(\scr{J}\) and must be the trivial
braid.  We have given a sketch of a proof of the following.

\begin{thm}\mylabel{IsoThmI} The group \(\widehat{BV}\) as defined
in this section and the group called \(\widehat{BV}\) as described
in Section \ref{GeomDescrSec} are isomorphic.  \end{thm}

\section{Distinguished subgroups}\mylabel{DisSubGpSec} 

In the group \(BV\), a single strand (corresponding to a single
Cantor set \(C\)) is split into a finite number \(n\) of strands
(corresponding to a cover of \(C\) by \(n\) intervals from
\tref{CoveringCollns}) which are then braided and recombined into
one strand.  This section picks out the appropriate subgroup of
\(\widehat{BV}\) and the corresponding subgroup of \(\widehat{V}\).

\subsection{Simple elements}\mylabel{SimplSect} We say that a hedge
\(H\) is {\itshape simple} if \(H(i)=1\) for all \(i>0\).  We say
that a forest \(F\) is {\itshape simple} if its corresponding hedge
\(c_F\) is simple.  Thus a simple forest \(F\) has at most one
non-trivial tree, and this non-trivial tree must be \(F_0\).  The
{\itshape type} of a simple forest \(F\) is the length of \(F\).
Note that the type of a simple forest \(F\) is also the number of
carets in \(F\) and is one less than the number of leaves of
\(F_0\).  Thus the type of the simple forest \(F\) is \(c_F(0)-1\).
The trivial forest is the only simple forest of type 0 and
\(\lambda_0\) is the only simple forest of type 1.

From Lemma \ref{FactorFacts}, we have that \(H(i)\le (HK)(i)\) for
all \(i\) for any hedges \(H\) and \(K\).  From this it follows that
if \(HK\) is simple, then \(H\) is simple.  Sending a forest \(F\)
to the corresponding hedge \(c_F\) is a homomorphism, so we get the
following.

\begin{lemma}\mylabel{LMSimple} If \(FG\) is simple for forests
\(F\) and \(G\), then \(F\) is simple.  \end{lemma}

If \(F\) is a simple forest of type \(k\), then the leaves of
\(F_0\) are numbered from \(0\) through \(k\).  We have that
\(F\lambda_i\) is simple if and only if \(i\le k\).  Inductively, we
get the following.

\begin{lemma}\mylabel{ForSimpleCond} A forest
\(F=\lambda_{i_1}\lambda_{i_2} \cdots \lambda_{i_k}\) is simple
of type \(k\) if
and only if \(i_j<j\) for all \(j\) with \(1\le j\le k\).
\end{lemma}

Recall that \(S_k\) (respectively, \(B_k\)) is the subgroup of
\(S_\infty\) (respectively, \(B_\infty\)) generated by \((\sigma_0,
\ldots, \sigma_{k-2})\).

If \(F\) is a forest and \(\beta\) is in \(S_\infty\) or
\(B_\infty\), then \(F\beta\) is {\itshape simple of type \(k\)} if
\(F\) is simple of type \(k\) and \(\beta\) is in \(S_{k+1}\) or
\(B_{k+1}\).  Intuitively, \(\beta\) permutes or braids only the
leaves of \(F_0\).

\begin{lemma}\mylabel{SimpleProdCond} Let \(F\) and \(G\) be in
\(\scr{F}\) and \(\beta\) and \(\gamma\) be from one of \(S_\infty\)
or \(B_\infty\) and assume that \(F\beta\) is simple of type \(k\)
and \(G=\lambda_{i_1}\cdots \lambda_{i_n}\).  Then the following are
equivalent.  {\claimenum \item \((F\beta)G\) is simple.  \item
\(F(\beta \cdot G)\) is simple.  \item \(i_j \le k+j-1\) whenever
\(1\le j\le n\).  \item \((F\beta)\lambda_{i_1} \cdots
\lambda_{i_j}\) is simple whenever \(1\le j\le n\).  \claimenumend}
Further, \((F\beta)(G\gamma)\) is simple if and only if
\((F\beta)G\) is simple and \(\gamma\) is in \(S_{k+n+1}\) or
\(B_{k+n+1}\).  \end{lemma}

\begin{proof} The definition of simple gives \((a)\Rightarrow (b)\).

We prove \((b)\Rightarrow (c)\Rightarrow (a)\) by induction on
\(n\).

If \(n=1\) and \(G=\lambda_i\), then \(\beta\cdot
\lambda_i=\lambda_{\beta(i)}\).  Since \(\beta\) can permute
non-trivially only the leaves of \(F_0\), we have that
\(F(\beta\cdot \lambda_i)\) is simple if and only if \(i\le k\).
Futher, when \(F(\beta\cdot \lambda_i)\) is simple, it is of type
\(k+1\).  

Now assume \(i\le k\).  If \(\beta=\beta'\sigma_j\) with \(\beta'\)
and \(\sigma_j\) in \(S_{k+1}\) or \(B_{k+1}\), then
\(\beta\lambda_i = \beta'\lambda_{\sigma_j(i)}
(\sigma_j)^{\lambda_i}\).  Since \(\sigma_j\) is in \(S_{k+1}\) or
\(B_{k+1}\), we have \(j\le k-1\) and \tref{ZapOpIII} gives us that
\((\sigma_j)^{\lambda_i}\) is in \(S_{k+2}\) or \(B_{k+2}\).  With
\(r=\sigma_j(i)\le k\) since \(j\le k-1\), we inductively get that
\((\beta')^{\lambda_r}\) is in \(S_{k+2}\) or \(B_{k+2}\) and thus
\((\beta)^{\lambda_i}\) is in \(S_{k+2}\) or \(B_{k+2}\).

We have proven \((b)\Rightarrow(c)\Rightarrow(a)\) in the case
\(n=1\).  The general case follows by induction and the equivalence
of \((a)\) and \((b)\) when \(n=1\).  The equivalence of \((d)\)
with the other statements is immediate.

The last claim follows from the equivalence of \((a)\) and \((b)\)
and the definitions.  \end{proof}

\begin{cor}\mylabel{SimpleProdEquiv} Let \(F\), \(F'\) and \(G\) be
in \(\scr{F}\) and let let \(\beta\), \(\beta'\) and \(\gamma\) be
from one of \(S_\infty\) or \(B_\infty\).  If \(F\beta\) and
\(F'\beta'\) are simple of the same type, then \((F\beta)(G\gamma)\)
is simple if and only if \((F'\beta')(G\gamma)\) is simple.
Further, if one (and thus both) of \((F\beta)(G\gamma)\) and
\((F'\beta')(G\gamma)\) is (are) simple, then they are of the same
type.  \end{cor}

\begin{proof} The first claim is a direct application of Lemma
\ref{SimpleProdCond} and the second follows from the definition of
type.  \end{proof}

\begin{lemma}\mylabel{SimpleCRM} Let \(F\) and \(G\) be in
\(\scr{F}\) and \(\beta\) and \(\gamma\) be both from either
\(S_\infty\) or \(B_\infty\).  If \(F\beta\) and \(G\gamma\)
are simple, then they have a simple common right multiple.
\end{lemma}

\begin{proof} This is easier than a reference to Lemma
\ref{PartMulProps} since the fact that \(S_\infty\) and \(B_\infty\)
are groups implies that any common right multiple of \(F\) and \(G\)
in \(\scr{F}\) is a common right multiple of \(F\beta\) and
\(G\gamma\).  We know that \(F\cup G\) is a common right multiple.
It is also clear that \(F\cup G\) is simple if both \(F\) and \(G\)
are simple.  \end{proof}

\subsection{Balanced, simple subgroups}\mylabel{BalSimplSect} The
groups \(\widehat{V}\) and \(\widehat{BV}\) are groups of fractions
and elements are represented by pairs of elements from a monoid.  We
pick out elements represented by certain pairs.

Let \((F\beta, F'\beta')\) be a pair of elements from
\(\FP\) or \(\FB\) with \(F\) and \(F'\) from \(\scr{F}\) and
\(\beta\) and \(\beta'\) from
\(S_\infty\) or \(B_\infty\) as appropriate.  We say that the pair
is {\itshape simple and balanced} if both entries in the pair are
simple, and if the two entries are of the same type.

Let \(V\) be the set of elements in \(\widehat{V}\) that have at
least one representative that is simple and balanced.  Let \(BV\) be
the set of elements in \(\widehat{BV}\) that have at least one
representative that is simple and balanced.  The point of Corollary
\ref{SimpleProdEquiv} and Lemma \ref{SimpleCRM} is the following.

\begin{thm}  Both \(V\) and \(BV\) are groups.
\end{thm}

\begin{proof} The inverse of a pair \((u, v)\) is \((v, u)\), so 
the groups of the statement are closed under inversion.  If
\((u,v)(w,z)\) is a product of pairs that must be calculated, then
we must find a common right multiple \(vp=wq\) of \(v\) and \(w\)
and get the product \((up, zq)\).  We know from Lemma
\ref{SimpleCRM} that \(vp=wq\) can be made simple and we know from
Corollary \ref{SimpleProdEquiv} that \(up\) and \(zq\) will be as
well.  Since type equals length of the forest part, since length of
forests is multiplicative, and since \(u\) and \(v\) share a type,
and \(w\) and \(z\) share a type, we get that the types of \(up\),
\(vp=wq\) and \(zq\) are the same.  Thus both groups of the
statement are closed under product.  \end{proof}

\begin{thm} The group \(BV\) from this section and the group \(BV\)
as described in Section \ref{GeomDescrSec} are isomorphic.
\end{thm}

\begin{proof}[Sketch of proof] The groups under discussion are
subgroups of the group \(\widehat{BV}\) realized as a group of
fractions and as described in Section \ref{GeomDescrSec}.  By
Theorem \ref{IsoThmI} the two versions of \(\widehat{BV}\) are
isomorphic.  The two versions of \(BV\) are the corresponding
subgroups.  \end{proof}

\providecommand{\bysame}{\leavevmode\hbox to3em{\hrulefill}\thinspace}

\noindent Department of Mathematical Sciences

\noindent State University of New York at Binghamton

\noindent Binghamton, NY 13902-6000

\noindent USA

\noindent email: matt@math.binghamton.edu

\end{document}